\documentclass[onefignu,onetabnumm]{siamart251216}

\usepackage{graphicx} 
\usepackage{algorithm}
\usepackage{algorithmic}

\usepackage[utf8]{inputenc}
\usepackage{graphicx}
\usepackage{bm}
\usepackage{comment}
\usepackage{psfrag}
\usepackage{latexsym,amsmath,amsfonts,amscd}
\usepackage{pifont}

\usepackage{changebar}
\usepackage{color}
\usepackage{bm}
\usepackage{tikz}
\usepackage{multirow}
\usepackage{xcolor}
\usepackage{hyperref}
\usepackage{setspace}
\usepackage{amssymb}
\usepackage{mathrsfs}
\usepackage{caption}
\usepackage{subcaption}

\usepackage[title]{appendix}
\usepackage{ulem}
\usepackage{tikz}

\usetikzlibrary{arrows,backgrounds,positioning}

\newtheorem{thm}{Theorem}[section]

\newtheorem{rem}[thm]{Remark}

\newcommand{\BU}{\mathbf{U}}
\newcommand{\BV}{\mathbf{V}}

\newcommand{\BP}{\mathbf{P}}
\newcommand{\BA}{\mathbf{A}}

\newcommand{\BR}{\mathbf{R}}

\newcommand{\BOmega}{\boldsymbol{\Omega}}

\newcommand{\bx}{\mathbf{x}}

\newcommand{\bw}{\mathbf{w}}
\newcommand{\bc}{\mathbf{c}}
\newcommand{\bb}{\mathbf{b}}

\newcommand{\bn}{\boldsymbol{n}}

\newcommand{\half}{\frac{1}{2}}
\newcommand{\br}{\boldsymbol{r}}
\newcommand{\bg}{\boldsymbol{g}}

\newcommand{\bphi}{\boldsymbol{\phi}}
\newcommand{\bpsi}{\boldsymbol{\psi}}
\newcommand{\BPsi}{\boldsymbol{\Psi}}

\newcommand{\bG}{\boldsymbol{G}}
\newcommand{\wbG}{\widetilde{\boldsymbol{G}}}

\newcommand{\BC}{{\bf{C}} }

\newcommand{\BM}{{\bf{M}} }
\newcommand{\BD}{{\bf{D}}}
\newcommand{\BX}{{\bf{X}}}

\newcommand{\BSigma}{{\mathbf{\Sigma}}}

\newcommand{\sk}{s^{(k)}}

\newcommand{\pzc}{\textcolor{black}}

\newcommand{\review}{\textcolor{black}}

\newcommand{\BS}{\mathbf{S}}

\newcommand{\algorithmicbreak}{\textbf{break}}
\newcommand{\BREAK}{\STATE \algorithmicbreak}

\title{Highly Efficient Rank-Adaptive Sweep-based SI-DSA for the Radiative Transfer Equation via Mild Space Augmentation }

\author{Wei Guo\footnote{Department of Mathematics and Statistics, Texas Tech University, Lubbock, TX, 70409, USA.},\quad   Zhichao Peng\footnote{Department of Mathematics, The Hong Kong University of Science and Technology, Clear Water Bay, Kowloon, Hong Kong, China. Corresponding author, E-mail: pengzhic@ust.hk.}}

\begin{document}

\maketitle

\begin{abstract}
Low-rank methods have emerged as a promising strategy for reducing the computational cost and storage of discrete-ordinates discretizations of the radiative transfer equation (RTE). However, most existing rank-adaptive approaches rely on rank-proportional space augmentation, which can negate efficiency gains when the effective solution rank becomes moderately large.
To overcome this limitation, we develop a rank-adaptive sweep-based source iteration with diffusion synthetic acceleration (SI–DSA) for the first-order steady-state RTE. The core of our method is a \review{sweep-based low-rank SI with an inner-loop iteration}  that performs efficient rank adaptation via mild space augmentation. In each inner iteration, the spatial basis is augmented with a small, rank-independent number of basis vectors without truncation, while a single truncation is performed only after the inner loop converges. Efficient rank adaptation is achieved through residual-based greedy random angular subsampling strategy together with incremental updates of projection operators, enabling non-intrusive reuse of existing transport-sweep implementations. In the outer iteration, a DSA preconditioner is applied to accelerate convergence.
Numerical experiments show that the proposed solver achieves accuracy and iteration counts comparable to those of full-rank SI–DSA while substantially reducing storage and runtime, even for challenging multiscale problems in which the effective rank reaches $30\text{–}45\%$ of the full rank.
\end{abstract}

\section{Introduction\label{sec:introduction}}
The radiative transfer equation (RTE) is a fundamental kinetic equation describing the transport and interaction of particles, such as photons and neutrons, as they propagate through and interact with a background medium. It has a wide range of applications in medical imaging, astrophysics and nuclear engineering. 
The unknown of RTE, commonly referred to as the angular flux or particle distribution depending on the application, is defined over a high-dimensional phase space comprising the spatial locations and angular directions of particle transport, \review{and, in many applications, additional frequency or energy variables.}  This high dimensionality renders deterministic discretizations computationally demanding in both memory and runtime. \review{In the present work, we restrict our attention to the space–angle formulation and do not consider additional energy or frequency variables.}

To mitigate the computational challenges, low-rank methods have been developed to approximate solutions of both steady-state and time-dependent RTE via matrix or tensor factorizations. Existing approaches  for time-dependent RTE mainly follow two paradigms: dynamical low-rank (DLR) methods \cite{koch2007dynamical}, which evolve the solution directly on a low-rank manifold by projecting the underlying dynamical system onto the tangent space; and step-and-truncation (SAT) strategies \cite{dektor2021rank,guo2022low}, which first evolve the solution in an enlarged space and then project it back onto a low-rank manifold via truncation or rounding. 
Instances of low-rank method for time dependent RTE include  
DLR methods for the first- or second-order formulations of RTE \cite{peng2020low,peng2023sweep,haut2026efficient}, DLR and SAT methods based on micro-macro decompositions \cite{einkemmer2021asymptotic,einkemmer2024asymptotic,frank2025asymptotic,sands2024high}, high-order low-order (HOLO) DLR method \cite{peng2021high}, and rank-adaptive predictor-corrector approaches \cite{hauck2023predictor}; these have been further analyzed in \cite{ding2021dynamical,yin2025towards}. \review{Under the SAT framework, tensor train decompositions \cite{oseledets2011tensor} have been exploited in low-rank methods for multigroup transport \cite{truong2023tensor,deshpande2026multigroup}.}
For steady state RTE, preconditioned rank-adaptive iterative solvers have been developed, including soft-thresholding based low-rank Richardson iteration with an exponential sum preconditioner \cite{bachmayr2024low} and a  low-rank source iteration with diffusion synthetic acceleration (SI-DSA) based on the second order reformulation of RTE \cite{guo2025inexact}. For reviews of low-rank methods for kinetic equations and high-dimensional PDEs, we refer to \cite{einkemmer2025review} and \cite{bachmayr2023low}, respectively.
 Despite these advances in low-rank methodologies, achieving substantial computational gains can be challenging, especially for practical problems whose solutions exhibit only moderately low-rank structure. The main reason is that some intermediate steps can be computationally expensive (e.g., scaling superlinearly in the solution rank) and may offset the efficiency gains of the low-rank approach. Consequently, the design of low-rank schemes that deliver genuine speedups beyond idealized test cases remains nontrivial.

In this work, we focus on the steady-state RTE, for which preconditioned iterative methods, such as source iteration (SI) and Krylov methods with diffusion synthetic acceleration (DSA) preconditioners, are well established \cite{Adams2002FastIM}. Although low-rank iterative solvers for the RTE remain comparatively underdeveloped, several pioneering efforts have emerged in recent years. Existing low-rank iterative methods can be broadly classified into two main categories. The first is based on a second-order reformulation of the RTE, within which low-rank methods have been developed under both rank-adaptive and DLR frameworks. Representative examples include the  preconditioned low-rank Richardson iteration in \cite{bachmayr2017iterative}, the rank-adaptive SI-DSA in \cite{guo2025inexact}, and the \review{$P_N$-like} DLR approach in \cite{haut2026efficient}. The second category directly targets the first-order formulation of the RTE and exploits transport sweeps, a fast matrix-free algorithm to solve the discrete streaming system. Existing methods along this direction include the projection-based fixed-rank DLR method in \cite{peng2023sweep} and the collocation-based  $S_N$-like DLR method in \cite{haut2026efficient}.
Despite these encouraging developments, several limitations remain:
\begin{enumerate}
    \item \underline{Ill-posedness of the second-order formulation:} When the total cross section is small or vanishes in a sub-region of the computational domain, the linear system resulting from the second-order reformulation of RTE becomes ill-posed or even singular, posing computational challenges to methods built upon it \cite{bachmayr2024low,guo2025inexact,haut2026efficient}.

    \item \emph{Overhead of  \review{space augmentation in rank-adaptation}:} \review{Current sweep-based low-rank methods in \cite{peng2023sweep,haut2026efficient} are fixed rank, however, the optimal rank required to balance accuracy and efficiency is typically unknown \textit{a priori}. Though rank adaptation mechanism such as basis update Galerkin (BUG) \cite{ceruti2022rank} could in principle be integrated into these approaches,} such  strategies often rely on $r$-dependent space augmentation, where $r$ denotes the solution rank. Typically, these methods enrich the mode space by increasing its dimension to a multiple of $r$ before truncation (e.g., $2r$ in BUG \cite{ceruti2022rank} or $4r$ for the SAT-based low-rank method in \cite{guo2025inexact} when solving  RTE \review{in the 2D X-Y geometry}). 
    However, such augmentation can be overly aggressive for challenging problems that require a relatively high solution rank to achieve a prescribed accuracy or satisfy the convergence criterion of the iterative solver. For example, the effective rank may reach approximately $40\%$ of the full rank \cite{guo2025inexact} in order to reduce the difference between the scalar fluxes of successive iterations below $10^{-5}$.   In these regimes, the computational overhead of $r$-dependent augmentation can offset performance gains, preventing genuine computational savings compared to efficient full-rank solvers.
\end{enumerate}

\review{To address these challenges in a focused setting, we consider steady-state linear radiative transfer for single-group neutron, photon, or proton transport,} and develop an efficient rank-adaptive low-rank SI-DSA
based on a mild augmentation procedure. Specifically, we design a novel rank-adaptation strategy that avoids rank-dependent, overly
aggressive basis enrichment. 
In our method, the enrichment is performed progressively, and the dimension of the resulting enriched mode space  is only slightly larger than $r$ prior to truncation, rather than a multiple thereof.
The key components of our approach are as follows.
\begin{enumerate}
    \item A novel iterative low-rank solver for the SI step is designed to achieve efficient rank-adaptation in an incremental manner. During each inner-loop iteration, a small, constant number of modes is added without immediate truncation; a single truncation is then performed upon convergence of the inner loop. This is realized through a novel residual-based greedy random subsampling of angular directions \review{balancing efficiency and sampling quality}.
    \item To further enhance efficiency of the low-rank SI method, the spatial low-rank basis and the projected operators are updated incrementally, thereby avoiding unnecessary global orthogonalizations and projections. 

    \item 
    Within SI-DSA, our low-rank SI serves as a seamless alternative to the full-rank SI, enabling the direct application of DSA and transport sweeps.

\end{enumerate}
\review{The primary contribution of this paper is not only to integrate rank adaptation into sweep-based low-rank iterative solvers for the steady-state first-order RTE, but also to design an efficient rank-adaptation strategy that avoids $O(r)$ space augmentation via randomized greedy angular sampling. This novel rank-adaptation approach via mild space augmentation enables computational speedup even for challenging problems with moderate-to-high solution ranks, e.g., $40\%$ of the full rank.}

We note that our method shares conceptual similarities with the collocation-based \review{$S_N$-like} DLR method recently proposed in the concurrent work \cite{haut2026efficient}. Both frameworks construct low-rank approximations by combining angular space sampling with projections in the physical space. However,  our approach differs fundamentally in the algorithmic integration of the angular sampling. The method in \cite{haut2026efficient} is designed for implicit time-marching: at each time step, angles are sampled based on the solution from the previous step, and a modified discrete ordinate system determined by these angles is subsequently solved via SI-DSA. Consequently, their sampling occurs externally to the SI-DSA iterative solver. In contrast, our method adaptively samples angles to update the solution rank internally within the inner-loop iterations of the low-rank SI. This integrated design enables our rank adaptation through mild, rather than aggressive, space augmentation.

Our paper is organized as follows. In Sec. \ref{sec:background}, we introduce the model problem and motivate the development of our method. In Sec. \ref{sec:lr-method}, we develop our  low-rank SI-DSA. In Sec. \ref{sec:numerical}, we numerically demonstrate the performance of our method through a series of benchmark tests. At last, we draw our conclusions in Sec. \ref{sec:conclusion}. 
\section{Background\label{sec:background}}
In this paper, we consider steady-state, single group, linear RTE with isotropic scattering and isotropic inflow boundary conditions on the computational domain $\Gamma_{\bx}$:
\begin{subequations}
\label{eq:rte}
    \begin{align}
    &\BOmega \cdot \nabla \psi(\bx,\BOmega) + \sigma_t(\bx) \psi(\bx,\BOmega) = \sigma_s(\bx) \phi(\bx) + G(\bx), \quad\sigma_t(\bx)=\sigma_s(\bx)+\sigma_a(\bx),\\
    &\phi(\bx) = \frac{1}{4\pi}\int_{\mathbb{S}^2} \psi(\bx,\BOmega) d\BOmega, \quad \bx\in\Gamma_{\bx},\;\BOmega\in\mathbb{S}^2,\\
    &\psi(\bx,\BOmega) = 0, \quad \bx\in \partial\Gamma_{\bx},\;\BOmega\cdot \bn(\bx)<0.\label{eq:boundary-condition}
    \end{align}
\end{subequations}
The RTE is posed in the $(\BOmega,\bx)$-phase space, where $\BOmega=(\Omega_x,\Omega_y,\Omega_z)\in\mathbb{S}^2$ and $\bx=(x,y,z)\in\Gamma_{\bx}$ denote the angular direction and spatial location, respectively. The angular flux at location $\bx$ in direction $\BOmega$ is denoted by $\psi(\bx,\BOmega)$, while the scalar flux is defined as $\phi(\bx)=\frac{1}{4\pi}\int_{\mathbb{S}^2} \psi(\bx,\BOmega) d\BOmega$. $G(\bx)$ represents an isotropic source. The scattering, absorption and total cross sections of the background medium are denoted by $\sigma_s(\bx)\geq 0$, $\sigma_a(\bx)\geq 0$ and $\sigma_t(\bx)=\sigma_s(\bx)+\sigma_a(\bx)$, respectively. The inflow boundary condition is given in \eqref{eq:boundary-condition}, where $\bn(\bx)$ denotes the outward normal direction of the computational domain at location $\bx\in\partial\Gamma_{\bx}$.

As the scattering cross section $\sigma_s\rightarrow\infty$, the solution of the RTE \eqref{eq:rte}  converges to its diffusion limit: $\psi(\bx,\BOmega)\rightarrow\phi(\bx)$, where $\phi(\bx)$ solves the diffusion equation
\begin{equation}
    -\nabla\cdot\left(\frac{1}{\sigma_s}\boldsymbol{\mathcal{D}}\nabla\phi\right)=-\sigma_a\phi+G,\;
    \boldsymbol{\mathcal{D}}=\frac{1}{4\pi}\textrm{diag}\left(\int_{\mathbb{S}^2}\Omega_x^2d\BOmega,\int_{\mathbb{S}^2}\Omega_y^2d\BOmega,\int_{\mathbb{S}^2}\Omega_z^2d\BOmega\right).
\end{equation}
\subsection{$S_N$ upwind discontinuous Galerkin discretization\label{sec:sn-dg}}
Following \cite{Adams2002FastIM,larsen2009advances}, we employ a discrete ordinates $(S_N)$ angular discretization \cite{pomraning2005equations} coupled with an upwind discontinuous Galerkin (DG) spatial discretization. As shown in \cite{adams2001discontinuous,guermond2010asymptotic},  this specific discretization is asymptotic preserving when utilizing piecewise linear or higher order polynomials. Consequently, the diffusion limit can be captured without requiring a refined mesh to resolve the particle mean free path. The details of our numerical scheme are presented below.

\textbf{$S_N$ angular discretization:} We solve RTE \eqref{eq:rte} at a set of quadrature points $\{\BOmega_j\}_{j=1}^{N_{\BOmega}}$ in the angular space and approximate the scalar flux through numerical quadrature. In particular, we seek $\psi(\bx,\BOmega_j)\approx\psi_j(\bx)$ for $j=1,\dots,N_{\BOmega}$,
\begin{subequations}
\label{eq:SN}
\begin{align}
&\BOmega_j\cdot\nabla\psi_j(\bx)+\sigma_t\psi(\bx)=\sigma_s\phi(\bx)+G(\bx),\quad\phi(\bx)=\sum_{j=1}^{N_{\BOmega}}\omega_j\psi_j(\bx),\;\bx\in\Gamma_{\bx},\\
&\psi_j(\bx)=g(\bx),\quad\bx\in\partial\Gamma_{\bx},\BOmega_j\cdot \bn(x)<0,
\end{align}
\end{subequations}
where $\{\omega_j\}_{j=1}^{N_{\BOmega}}$ denote the normalized quadrature weights associated with $\{\BOmega_j\}_{j=1}^{N_{\BOmega}}$ satisfying $\sum_{j=1}^{N_{\BOmega}}\omega_j=1$. In this paper, we employ the Chebyshev–Legendre quadrature; see Appx. \ref{sec:angular-discretization} for details.

\textbf{Upwind DG method to solve the $S_N$ system \eqref{eq:SN}:} For simplicity, we focus on a 2D rectangular spatial domain $\Gamma_x=[x_L,x_R]\times[y_B,y_T]$ partitioned by a Cartesian mesh $\mathcal{T}_h=\{\mathcal{T}_{a,b}=[x_{a-\half},x_{a+\half}]\times[y_{b-\half},y_{b+\half}],\;1\leq a\leq N_x,\;1\leq b\leq N_y\}$ with $x_{a-\half}=x_L+(a-\frac{1}{2})\frac{x_R-x_L}{N_x}$ and $y_{b-\half}=y_B+(b-\half)\frac{y_T-y_B}{N_y}$. Specifically, $x_{\half}=x_L$, $x_{N_x+\half}=x_R$, $y_\half=y_B$ and $y_{N_y+\half}=y_T$. 
For $j=1,2,\dots,N_{\BOmega}$, we seek the numerical solution $\psi(\bx,\BOmega_j)\approx\psi_{h,j}\in U_h^K(\mathcal{T}_h)$ in the discrete space
\begin{subequations}
    \begin{align}
    &U_h^K(\mathcal{T}_h):=\{u(\bx): u(\bx)|_{\mathcal{T}_{i}}\in Q^K(\mathcal{T}_{a,b}),1\leq a \leq N_x,1\leq b\leq N_y\},\\
    &Q^K(\mathcal{T}_{a,b})=\{x^py^q,0\leq p,q\leq K,(x,y)\in \mathcal{T}_{a,b}\}.\label{eq:discrete_space}
    \end{align}
\end{subequations}
Here, $Q^K(\mathcal{T}_{a,b})$ denotes the space of tensor-product polynomials of degree at most $K$ in each coordinate on the cell $\mathcal{T}_{a,b}$.
More specifically, the functions $\psi_{h,j}$ satisfy $\forall \eta_h(\bx)\in U_h^K(\mathcal{T}_h),$
\begin{align*}
&\sum_{a=1}^{N_x}\sum_{b=1}^{N_y}\Bigg(-\int_{\mathcal{T}_{a,b}} \Big(\BOmega_j\cdot\nabla\eta_h\Big) \psi_{h,j} d\bx+\int_{\mathcal{T}_{a,b}}\sigma_t \psi_{h,j}\eta_hd\bx
\notag\\
&+\int_{y_{b-\half}}^{y_{b+\half}}\Omega_{j,x}\widehat{\psi}_{h,j}(x_{a+\half},y)\eta_h(x_{a+\half}^-,y)dy
-\int_{y_{b-\half}}^{y_{b+\half}}\Omega_{j,x}\widehat{\psi}_{h,j}(x_{a+\half},y)\eta_h(x_{a+\half}^+,y)dy
\notag\\
&+\int_{x_{a-\half}}^{y_{a+\half}}\Omega_{j,y}\widehat{\psi}_{h,j}(x,y_{b+\half})\eta_h(x,y_{b+\half}^-)dx
-\int_{x_{a-\half}}^{x_{a+\half}}\Omega_{j,y}\widehat{\psi}_{h,j}(x,y_{b-\half})\eta_h(x,y_{b-\half}^+)dx
\Bigg)
\notag\\
&= \sum_{a=1}^{N_x}\sum_{b=1}^{N_y}\left(\int_{\mathcal{T}_{a,b}}\sigma_s(\bx) \phi_h(\bx)\eta_h(\bx) d\bx
 + \int_{\mathcal{T}_{a,b}}G(\bx)\eta_h(\bx) d\bx\right),
\end{align*}
where $\phi_h=\sum_{j=1}^{N_{\BOmega}}\omega_j\psi_{h,j}$, $h(\alpha^{\pm},y)=\lim_{x\rightarrow \alpha^{\pm}}h(\alpha,y)$, $h(x,\beta^{\pm})=\lim_{y\rightarrow \beta^{\pm}}h(x,y)$  and the upwind flux $\widehat{\psi}_{h,j}$ is  defined as
\begin{subequations}
    \begin{align}
        &\Omega_{j,x}\widehat{\psi}_{h,j}(x_{a-\half},y)=\Omega_{j,x}\psi_{h,j}(x_{a-\half}^\mp,y),\quad\text{if }\pm\Omega_{j,x}\geq 0,\; 1\leq a\leq N_x+1,\\
        &\Omega_{j,y}\widehat{\psi}_{h,j}(x,y_{b-\half})=\Omega_{j,y}\psi_{h,j}(x,y_{b-\half}^\mp),\quad\text{if }\pm\Omega_{j,y}\geq 0,\;1\leq b\leq N_y+1.
    \end{align}
\end{subequations}
In addition, the boundary fluxes are determined by $\psi_{h,j}(x_\half^-,y)=g(x_L,y)$ for $\Omega_{j,x}\geq 0$, $\psi_{h,j}(x_{N_x+\half}^+,y)=g(x_R,y)$ for $\Omega_{j,x}\leq 0$, $\psi_{h,j}(x,y_\half^-)=g(x,y_B)$ for $\Omega_{j,y}\geq 0$ and $\psi_{h,j}(x,y_{N_y+\half}^+)=g(x,y_T)$ for $\Omega_{j,y}\leq 0$. 

\textbf{Matrix-vector formulation of the fully discrete system.} Let $\{\eta_k\}_{i=1}^{N_{\bx}}$ with $N_{\bx}=N_xN_y(K+1)^2$ be an orthonormal basis for $U_h^K$, and assume that $\eta_i$ is locally supported in the element $\mathcal{T}_{a_i,b_i}$, i.e., $\eta_i(x,y)=0$ if $(x,y)\not\in\mathcal{T}_{a_i,b_i}$. 
The numerical solution $\psi_{h,j}$ can be expanded as $\psi_{h,j}=\sum_{i=1}^{N_{\bx}}\psi_{ij}\eta_i$. Let $\bpsi_j=(\psi_{1j},\dots,\psi_{N_{\bx}j})^T\in\mathbb{R}^{N_{\bx}}$, which satisfies the following linear system: 
\begin{subequations}  
\begin{align}&(\BD_j+\BSigma_t)\bpsi_j=\BSigma_s\bphi+\bG+\bg_j^{\textrm{bc}}=\BSigma_s\bphi+\wbG_j,\quad \bphi=\sum_{j=1}^{N_{\BOmega}}\omega_j\bpsi_j, \label{eq:mat-vec-a}\\
\text{where}\;&\BD_j=\begin{cases}
\Omega_{j,x}\BD_x^-+\Omega_{j,y}\BD_y^-,\;\text{if }\Omega_{j,x}\geq 0,\Omega_{j,y}\geq 0,\\
\Omega_{j,x}\BD_x^-+\Omega_{j,y}\BD_y^+,\;\text{if }\Omega_{j,x}\geq0,\Omega_{j,y}<0,\\
\Omega_{j,x}\BD_x^++\Omega_{j,y}\BD_y^-,\;\text{if }\Omega_{j,x}<0,\Omega_{j,y}\geq0,\\
\Omega_{j,x}\BD_x^++\Omega_{j,y}\BD_y^+,\;\text{otherwise}.
\end{cases} \label{eq:upwind-operator}
\end{align}
\label{eq:mat-vec}
\end{subequations}
Here, $\BD_j\in\mathbb{R}^{N_{\bx}\times N_{\bx}}$ denotes the discrete transport operator associated with $\BOmega_j\cdot\nabla$, 
$\BSigma_s,\BSigma_t\in\mathbb{R}^{N_{\bx}\times N_{\bx}}$ denote the discrete scattering and total cross-section operators, and $\bG,\,\bg_j^{\textrm{bc}}\in\mathbb{R}^{N_{\bx}}$ correspond to the discretizations of the source term and the boundary conditions, respectively. 
For brevity, the precise definitions of $\BD_x^{\pm}, \BD_y^{\pm}, \BSigma_t, \BSigma_s, \bG,$ and $\bg_j^{\mathrm{bc}}$ are provided in Appx.~\ref{app:mat-vec}.

\subsection{Transport sweep and source iteration with diffusion synthetic acceleration (SI-DSA)}\label{sec:fr_sidsa}
SI-DSA \cite{Adams2002FastIM,larsen2009advances} is a popular preconditioned  iterative method to solve \eqref{eq:mat-vec}. A crucial component of the SI framework is the transport sweep \cite{Adams2002FastIM,larsen2009advances}, a highly efficient, matrix-free algorithm designed to invert the discrete streaming operator $\BD_j+\BSigma_t$. Below, we briefly review the fundamentals of SI-DSA and the associated matrix-free transport sweeps.

\textbf{SI-DSA.}  Given an initial guess for the scalar flux $\bphi^{(0)}$, the $n$-th iteration of SI-DSA follows a two-step procedure.
\begin{enumerate}
    \item \textbf{The SI step updating the angular flux}: We update the angular flux by solving the transport equations for each discrete ordinate independently, using the scalar flux from the previous iteration:
    \begin{equation}(\BD_j+\BSigma_t)\bpsi_j^{(n)}=\BSigma_s\bphi^{(n-1)}+\wbG_j,\label{eq:SI}
    \end{equation}
    and then update the scalar flux as $\bphi^{(n,*)}=\sum_{j=1}^{N_{\BOmega}}\omega_j\bpsi_j^{(n)}$.
    \item \textbf{The DSA step accelerating convergence}. DSA introduces a correction to the scalar flux by solving:
    \begin{equation}
    -\BD_{\textrm{diff.}}\delta\bphi^{(n)}=-\BSigma_a\delta\bphi^{(n)}+\BSigma_s(\bphi^{(n,*)}-\bphi^{(n)}),
    \label{eq:DSA}
    \end{equation}
    and then corrects scalar flux as $\bphi^{(n)}=\bphi^{(n,*)}+\delta\bphi^{(n)}$.
    Equation \eqref{eq:DSA} represents the discretization of the diffusion limit for the error equation
    \begin{equation}
        -\frac{1}{3}\nabla\cdot\left(\frac{1}{\sigma_s}\nabla\delta\phi^{(n)}\right)=-\sigma_a\delta\bphi^{(n)}+\sigma_s(\phi^{(n,*)}-\phi^{(n)}). 
        \label{eq:diffusion-correction}
    \end{equation}
We emphasize that, as discussed in \cite{Adams2002FastIM,larsen2009advances,palii2020convergent}, a consistent discretization of \eqref{eq:diffusion-correction} is required to ensure effective acceleration. Specifically, we employ a fully consistent discretization following \cite{adams2001discontinuous}, with detailed definitions provided in Appx.~B of \cite{peng2024reduced}. We also note that, as shown in \cite{larsen2009advances}, SI may converge arbitrarily slowly without DSA.

\end{enumerate}
The convergence criteria of SI-DSA is  chosen as $\|\bphi^{(n)}-\bphi^{(n-1)}\|_\infty\leq \epsilon_{\textrm{SI-SA}}$, where $\epsilon_{\textrm{SI-SA}}$ is a prescribed threshold.

\textbf{Transport sweeps.} Transport sweeps leverage the following property to solve \eqref{eq:SI} very efficiently:  with upwind discretization, the discrete streaming operator,  $\BD_j+\BSigma_t$, becomes block lower triangular when elements are ordered along the upwind direction corresponding to $\BOmega_j$. Hence,  the SI equation \eqref{eq:SI} can be solved very efficiently in a matrix-free manner through a single block Gauss–Seidel iteration sweeping elements along the upwind direction.

To better illustrate the transport sweep mechanism, we describe the block structure and sweeping procedure for the RTE in 1D slab geometry:
\begin{equation}
    v\partial_x \psi +\sigma_t \psi = \sigma_s\phi +G,\; x\in[x_L,x_R], v\in[-1,1]
\end{equation}
with inflow boundary conditions.
We partition the computational domain as 
$$[x_L,x_R]=\cup_{i=1}^{N_x}I_a,\;I_a=[x_L+(a-1)\Delta x,x_R+a\Delta x],\;\Delta x=(x_R-x_L)/N_x$$ 
and  apply $S_N$ upwind DG discretization using piecewise polynomials of degree at most $K$. The resulting linear system retains the form \eqref{eq:mat-vec-a}. Define $\bpsi^{(n)}_{j,a}\in\mathbb{R}^{K+1}$ as the vector of DOFs of $\bpsi_j^{(n)}$  associated with local basis functions in $I_a$. By ordering the elements along the upwind direction corresponding to $v_j$, namely from $1$ to $N_x$ for $v_j>0$ and from $N_x$ to $1$ for $v_j<0$, the SI equation takes a block lower triangular form:
\begin{subequations}
    \begin{align}
     &\left(\begin{matrix}
     \BA^{[j]}_{1,1} & \\
     \BA^{[j]}_{2,1} & \BA^{(v_j)}_{2,2} & \\
     & \ddots  & \ddots  & \\
             &  & \BA^{[j]}_{N_x,N_x-1} & \BA^{[j]}_{N_x,N_x}
     \end{matrix}\right)   
     \left(
     \begin{matrix}
         \bpsi^{(n)}_{j,1}\\
         \bpsi^{(n)}_{j,2}\\
         \vdots\\
         \bpsi^{(n)}_{j,N_x}
     \end{matrix}
     \right)
     =     \left(
     \begin{matrix}
         \bb^{(n)}_{j,1}\\
         \bb^{(n)}_{j,2}\\
         \vdots\\
         \bb^{(n)}_{j,N_x}
     \end{matrix}
     \right),\;\text{with }v_j\geq 0,\\
    &\left(\begin{matrix}
     \BA^{[j]}_{N_x,N_x} & \\
     \BA^{[j]}_{N_x-1,N_x} & \BA^{[j]}_{N_x-1,N_x-1} & \\
     & \ddots  & \ddots  & \\
             &  & \BA^{[j]}_{1,2} & \BA^{[j]}_{1,1}
     \end{matrix}\right)   
     \left(
     \begin{matrix}
         \bpsi^{(n)}_{j,N_x}\\
         \bpsi^{(n)}_{j,N_x-1}\\
         \vdots\\
         \bpsi^{(n)}_{j,1}
     \end{matrix}
     \right)
     =     \left(
     \begin{matrix}
         \bb^{(n)}_{j,N_x}\\
         \bb^{(n)}_{j,N_x-1}\\
         \vdots\\
         \bb^{(n)}_{j,1}
     \end{matrix}
     \right),\; \text{with }v_j<0.
    \end{align}
\end{subequations}
As a result, $\bpsi_j^{(n)}$ can be computed through a single transport sweep, which is equivalent to a single block Gauss-Seidel iteration:
\begin{itemize}
   \item For $v_j\geq0$, sweep from $a=1$ to $a=N_x$ and solve   $\BA_{a,a}^{[j]}\bpsi_{j,a}^{(n)}=\bb^{(n)}_{j,a}-\BA_{a,a-1}^{[j]}\bpsi_{j,a-1}$ with $\bpsi_{j,0}^{(n)}=\mathbf{0}$.
  \item  For $v_j<0$, sweep from $a=N_x$ to $a=1$ and solve $\BA_{a,a}^{[j]}\bpsi_{j,a}^{(n)}=\bb^{(n)}_{j,a}-\BA_{a,a+1}^{[j]}\bpsi_{j,a+1}$  with $\bpsi_{j,N_x+1}^{(n)}=\mathbf{0}$.
\end{itemize}
 This sweeping procedure can be straightforwardly extended to a 2D rectangular computational domain partitioned by a Cartesian mesh.

\subsection{Low-rank solvers for steady-state RTE and challenges}
To mitigate the substantial memory and computational costs arising from the high dimensionality of the phase space, low-rank methods approximate the discrete solution using low-rank matrix or tensor decompositions. For example,
\begin{equation}
\begin{bmatrix}
\bpsi(\cdot,\BOmega_1) & \bpsi(\cdot,\BOmega_2) & \dots & \bpsi(\cdot,\BOmega_{N_{\BOmega}}
)
\end{bmatrix}
\approx 
\boldsymbol{\Psi} =
\begin{bmatrix}
\bpsi_1 & \bpsi_2 & \dots & \bpsi_{N_{\BOmega}}
\end{bmatrix}
\approx
\BX \BS \BV^T,
\label{eq:lr-representation}
\end{equation}
where $\boldsymbol{\Psi}\in\mathbb{R}^{N_{\bx}\times N_{\BOmega}}$, $\BX\in\mathbb{R}^{N_{\bx}\times r}$, $\BS\in\mathbb{R}^{r\times r}$, and $\BV\in\mathbb{R}^{N_{\BOmega}\times r}$. 
The low-rank factors $\BX$, $\BS$, and $\BV$ can be computed using low-rank solvers based on DLR methods \cite{koch2007dynamical} or rank-adaptive approaches \cite{bachmayr2023low,guo2022low}. 
When the solution exhibits strong low-rank structure, i.e., $r\ll\min(N_{\bx},N_{\BOmega})$, these methods often achieve substantial savings in both \review{storage} and computational cost.

In this subsection, we first briefly review existing low-rank iterative solvers for the RTE. We then motivate the present work by highlighting the necessity of, and challenges in, designing highly efficient rank-adaptive low-rank iterative solvers.

\subsubsection{Existing low-rank iterative solvers} 
Low-rank iterative solvers for RTE have been developed 
to solve linear systems arising from steady-state RTE or implicit time marching of time-dependent RTE. Existing methods can be broadly divided into the following two categories.

\textbf{Low-rank solvers for the second-order even-parity reformulation of RTE.} Methods along this line include \cite{bachmayr2017iterative,guo2025inexact,haut2026efficient}. Rank-adaptive truncated summation has been utilized in low-rank Richardson iteration with an exponential sum preconditioner \cite{bachmayr2017iterative} and in low-rank SI--DSA methods for the second-order RTE \cite{guo2025inexact}. A fixed-rank $P_N$-like DLR method for the second-order RTE was introduced in \cite{haut2026efficient}.

\textbf{Sweep-based fixed-rank DLR methods for the first-order formulation of the RTE.}
An alternative strategy is to directly design a low-rank solver for the first-order RTE based on SI and transport sweeps \cite{peng2023sweep,haut2026efficient}. 
Peng \textit{et al.}~\cite{peng2023sweep} develop a fixed-rank DLR approach by grouping angular fluxes corresponding to directions within the same octant. This strategy requires modifications to the transport sweeps to solve projected low-rank systems for each angular group. 
To avoid modifying transport sweeps, Haut \textit{et al.}~\cite{haut2026efficient} propose a collocation-based \review{$S_N$-like} DLR method \cite{dektor2025interpolatory} that constructs low-rank solutions via angular sampling using the discrete empirical interpolation method (DEIM) \cite{sorensen2016deim} together with Galerkin projections in the physical space.

\subsubsection{Necessity and challenges to design rank-adaptive sweep-based solver\label{sec:motivation}}
The efficiency and accuracy of low-rank methods are governed by the effective solution rank, which is generally unknown \textit{a priori}. Therefore, adaptively determining the solution rank is essential for achieving genuine computational savings \review{without scarificing accuracy}. Standard rank-adaptive strategies enrich the mode space by expanding its dimension from $r$ to a multiple of $r$, followed by truncation \cite{bachmayr2023low,ceruti2022rank}. For instance, the rank-adaptive  basis-update-Galerkin (BUG) DLR method enlarges the space to $2r$ \cite{ceruti2022rank}, while the SAT method in \cite{guo2022low} for the Vlasov equation expands it to $3r$ in 1D1V. In addition, the truncation step typically involves operations that scale superlinearly in $r$, such as QR or SVD factorizations. When the effective rank required to achieve the desired accuracy or satisfy the convergence criteria of the iterative solver becomes moderately large, such rank-dependent augmentation can introduce substantial computational overhead. In these regimes, the aggressive space augmentation may offset or even negate the potential savings of low-rank methods compared with efficient full-rank solvers.
Indeed, for several challenging benchmark problems considered in this work, the effective ranks are relatively high. To reach the SI-DSA stopping criterion $\|\bphi^{(n)}-\bphi^{(n-1)}\|_\infty \le 10^{-6}$, the effective rank is slightly larger than $45\%$ of the full rank $\min\{N_{\bx},N_{\BOmega}\}$ for the lattice problem (Sec.~\ref{sec:lattice}), $36\%$ for the variable scattering case (Sec.~\ref{sec:variable-scattering}), and $40\%$ for the pin-cell case (Sec.~\ref{sec:pin-cell}). In such cases, even augmenting the mode space from dimension $r$ to $2r$ may render the low-rank solver slower than its full-rank counterpart.

Therefore, we aim to develop a rank-adaptive sweep-based low-rank iterative solver for the RTE that requires only mild space augmentation, thereby achieving practical reductions in both storage and computational cost across a broad range of problems, even when the effective rank reaches $30$–$40\%$ of the full rank.

\section{Rank-adaptive sweep-based SI-DSA via mild space augmentation\label{sec:lr-method}}

In this section, we first give an overview of our low-rank SI algorithm in Sec. \ref{sec:alg-overview}, then delve into details of crucial components, including the residual-based greedy random subsampling strategy for angular directions (Sec.\ref{sec:sampling}), the incremental update of low-rank bases and projected operators (Sec.\ref{sec:lr-update}), and the stopping criteria for the inner loop of the low-rank SI (Sec.\ref{sec:stop}). We summarize our contributions and highlight the key differences between the proposed method and existing sweep-based low-rank approaches \cite{peng2023sweep,haut2026efficient} in Sec.\ref{sec:contribution}. Finally, we outline possible extensions of the method to time-dependent problems in Sec.\ref{sec:extension}.

\subsection{Algorithm overview of iterative low-rank SI\label{sec:alg-overview}}
In this work, we employ a matrix-based low-rank representation for  the angular flux:
\begin{equation}
    \boldsymbol{\Psi}=
 \begin{bmatrix}
\bpsi_1,&\bpsi_2,&\dots& \bpsi_{N_{\BOmega}}   
\end{bmatrix} \approx
\BX\BS\BV^T,
\end{equation}
where $\BX\in\mathbb{R}^{N_{\bx}\times r}$, $\BS\in\mathbb{R}^{r\times r}$ and $\BV\in\mathbb{R}^{N_{\BOmega}\times r}$. 

Our low-rank SI seeks an approximate solution of the SI equation \eqref{eq:SI} in low-rank form, \(\BPsi^{(n)} \approx \BX^{(n)} \BS^{(n)} (\BV^{(n)})^T\), via an iterative procedure. In each inner-loop iteration, both the spatial basis and the rank are updated by augmenting the mode space with a small, fixed number of solutions to equation \eqref{eq:SI} corresponding to selected angular directions. A residual-based greedy random subsampling strategy is then employed to select new angular directions for the subsequent iteration.
Upon the convergence of the inner-loop, we construct the low-rank angular basis $\BV^{(n)}$ and the coefficient core $\BS^{(n)}$.     Given $\BX^{(n,k-1)}\in\mathbb{R}^{N_{\bx}\times r_{k-1}}$, the $k$-th inner-loop iteration obtains $\BX^{(n,k)}$ as follows: 
\begin{enumerate}
    \item \textbf{Transport sweeps:} Given $p$ newly sampled angles, $\{\BOmega_{s^{(k)}_1},\dots,\BOmega_{s^{(k)}_p}\}$, we apply transport sweeps to solve 
    \begin{equation}
        (\BD_j+\BSigma_t)\bpsi_j^{(n)}=\BSigma_s\bphi^{(n-1)}+\wbG_j,\quad j=s_1^{(k)},\dots,s_p^{(k)}.
        \label{eq:sampled-sweep}
    \end{equation}
    In our numerical experiments, we set $p\leq 3$.
    \item \textbf{Incremental update of low-rank spatial basis and operators.} We update \(\BX^{(n,k)}\) using the \review{modified Gram–Schmidt (MGS) with double projection and adaptive reorthogonalization following \cite{daniel1976reorthogonalization}}, as described in Alg.~\ref{alg:mgs2-a}, to compute its QR factorization incrementally.
    \begin{equation}
    \BX^{(n,k)}\BR^{(n,k)}=
    \begin{bmatrix}
        \BX^{(n,k-1)} & \bpsi^{(n)}_{\sk_1} & \dots & \bpsi^{(n)}_{\sk_p}
    \end{bmatrix}, 
    \label{eq:qr}
    \end{equation}
    where $\BX^{(n,k)}\in\mathbb{R}^{N_{\bx}\times r_k}$ has orthonormal columns, $\BR\in\mathbb{R}^{r_k\times r_k}$ is upper triangular, and the updated rank satisfies $r_k=r_{k-1}+p$. \review{Particularly, by adaptively applying additional MGS projections to mitigate the loss of orthogonality \cite{daniel1976reorthogonalization},}  $\BX^{(n,k)}$ and $\BX^{(n,k-1)}$ satisfy the hierarchical relation $\BX^{(n,k)}=\begin{bmatrix} \BX^{(n,k-1)} & \BP^{(n,k)}\end{bmatrix}$. Leveraging this hierarchical relation, the projection operators used in the sampling step can be updated efficiently in an incremental manner. See Sec. \ref{sec:lr-update}. 
    \item \textbf{Residual-based greedy random subsampling of angular directions.}  We randomly sample $q>p$ previously unsampled angles to form the candidate set $\mathcal{S}^{(k)}=\{\BOmega_{t^{(k)}_1},\dots,\BOmega_{t^{(k)}_q}\}$.  We obtain reduced solutions for angles in $\mathcal{S}^{(k)}$ through Galerkin projections. Based on the corresponding PDE residuals, we greedily select the $p$ angles in $\mathcal{S}^{(k)}$ with the largest residuals $\{\BOmega_{s^{(k+1)}_1},\dots,\BOmega_{s^{(k+1)}_p}\}$. See Sec. \ref{sec:sampling}.
    \item \textbf{Inner-loop stopping criteria.} If the proposed  stopping criteria in Sec. \ref{sec:stop} are satisfied, we compute $\phi^{(n,*)}$ via numerical integration, and update the low-rank bases and coefficient core through projections and SVD truncation. 
\end{enumerate}
Main features of the above low-rank SI are summarized as follows.
\begin{enumerate}
\item \textbf{Sampling in the angular space and Galerkin projection in the physical space.}
We combine angular sampling with physical-space Galerkin projections to construct low-rank solutions, similar to \cite{haut2026efficient}. This hybrid strategy enables the reuse of existing transport-sweep implementations in a minimally invasive manner; see Sec.~\ref{sec:contribution} for details.

\item\textbf{Mild space augmentation to achieve rank adaptivity.} In each inner-loop iteration of our low-rank SI,  the mode space is mildly enlarged by adding a small rank-independent number, namely $p$, of snapshots without immediate truncation. A single truncation step is performed only after the inner loop has converged to compress the low-rank representation. As a result, the dimension of the augmented mode space before one-time truncation is only slightly larger than $r$, in contrast to most existing methods that rely on overly aggressive $O(r)$ space augmentation.

\item\textbf{Delayed construction of $\BS^{(n)}$ and $\BV^{(n)}$ to enable efficient subsampling.} We observe that, with a properly designed stopping criterion, it is unnecessary to explicitly construct $\BV^{(n)}$ and $\BS^{(n)}$ at every inner-loop
iteration of the low-rank SI. This design choice enables an efficient
residual-based greedy random subsampling strategy for rank adaptation. The
proposed subsampling approach reduces the cost of solving projected systems
during angular sampling from $O(N_{\BOmega} r^3)$ to $O(q r^3)$ with
$q \ll N_{\BOmega}$. 

\item \textbf{Incremental update of spatial basis and projected operators.} Since no truncation is performed before the convergence of the inner-loop
iterations, the adaptively expanded mode space used to construct the low-rank spatial basis possesses a hierarchical structure. Leveraging this structure, we perform efficient incremental updates of both the low-rank spatial basis and the associated Galerkin projected operators, thereby further improving computational efficiency.

\item \review{\textbf{Reuse of existing one-angle sweep implementations}. Since no modification is introduced to the SI equation for sampled angles \eqref{eq:sampled-sweep}, existing sweep implementations can be directly applied.}
\end{enumerate}

With the proposed rank-adaptive low-rank SI, we are now ready to introduce the
corresponding low-rank SI-DSA scheme. Each iteration still follows the same two-step procedure as in the full-rank setting: we first apply
the low-rank SI to compute a low-rank approximation of the SI equation
\eqref{eq:SI}, and then correct the scalar flux to accelerate convergence by solving the DSA equation
\eqref{eq:DSA}. The outer-loop stopping criterion is still chosen as $\|\bphi^{(n)}-\bphi^{(n-1)}\|_\infty$ to be sufficiently small. Details of the low-rank SI-DSA
 are summarized in Alg.~\ref{alg:lr-si-dsa}.

\begin{algorithm}[ht]
\caption{Low-rank SI-DSA.\label{alg:lr-si-dsa}}
\begin{algorithmic}[1]
\STATE{\textbf{Input}: Initial scalar flux $\bphi^{(0)}$, tolerance $\epsilon_{\textrm{Diff}}$, and maximum number of iterations, $N_{\textrm{iter}}$.}
\STATE{\textbf{Initialization}: Set $n=1$.}
\WHILE{$k\leq N_{\textrm{Iter}}$}
\STATE{Apply the \textbf{low-rank SI} following Sec. \ref{sec:alg-overview} to solve \eqref{eq:SI} and obtain $\bphi^{(n,*)}$. Optionally return a rank-$r_n$ factorization of the angular flux,
$\BX^{(n)} \BS^{(n)} (\BV^{(n)})^T$.
}

\IF{$\|\bphi^{(n,*)}-\bphi^{(n-1)}\|_{\infty}\leq\epsilon_{\textrm{SI-SA}}$}
    \BREAK
\ENDIF
\STATE{Compute the \textbf{DSA correction} $\delta\bphi^{(n)}$ by solving a diffusion correction equation, then update the scalar flux as $\bphi^{(n+1)}=\bphi^{(n,*)}+\delta\bphi^{(n)}$.
\STATE Set $n \gets n+1$.
}
\ENDWHILE
\STATE{\textbf{Output:} Low-rank angular flux $\BX^{(n)} \BS^{(n)} (\BV^{(n)})^T$ and scalar flux $\bphi^{(n)}$.}
\end{algorithmic}
\end{algorithm}
\subsection{Greedy random subsampling of angles\label{sec:sampling}}
We present the greedy random subsampling strategy to select angles for rank adaptation in the subsequent iteration. As outlined in Sec. \ref{sec:alg-overview}, instead of considering all angular directions, we first randomly sample $q$ previously unsampeld angles to form the candidate set $\mathcal{S}^{(k)}=\{\BOmega_{t^{(k)}_1},\dots,\BOmega_{t^{(k)}_q}\}$, then compute reduced solutions corresponding to angles in $\mathcal{S}^{(k)}$ through Galerkin projections. Based on the associated residuals, we greedily pick $p$ angular directions with the largest residuals for augmentation in the next iteration.

\textbf{Galerkin projection.}
Given the updated spatial basis $\BX^{(n,k)}$, for an angular direction $\BOmega_j\in\mathcal{S}^{(k)}$ we compute the reduced solution via a Galerkin projection
\begin{equation}
\BA_{r_k}^{(n,k)}\bc_j^{(n,k)}
\triangleq
(\BX^{(n,k)})^T(\BD_j+\BSigma_t)\BX^{(n,k)}\bc_j^{(n,k)}
=
(\BX^{(n,k)})^T(\BSigma_s\bphi^{(n-1)}+\wbG_j),
\label{eq:galerkin-projection}
\end{equation}
where $\bc_j^{(n,k)}\in\mathbb{R}^{r_k}$ and $r_k=pk$ denotes the rank of the spatial basis $\BX^{(n,k)}$. The projected approximation of the angular flux corresponding to these directions is given by $\bpsi_j^{(n,k)}\approx \BX^{(n,k)}\bc_j^{(n,k)}$.

\textbf{Residual-based greedy sampling.}
After computing the projected solutions, we evaluate the PDE residual of the SI equation \eqref{eq:SI} for each candidate angular direction in $\mathcal{S}^{(k)}$:
\begin{equation}
\br^{(k)}_j
= \BSigma_s \bphi^{(k)} + \wbG_j
- (\BD_j + \BSigma_t)\BX^{(n,k)} \bc_j^{(n,k)},
\qquad \Omega_j \in \mathcal{S}^{(k)}.
\end{equation}
We then sort the indices of the angles in $\mathcal{S}^{(k)}$ in descending order of their $\ell_2$-norms:
\[
\|\br^{(k)}_{s^{(k+1)}_1}\|_2
\ge
\cdots
\ge
\|\br^{(k)}_{s^{(k+1)}_q}\|_2.
\]
The first $p$ angles corresponding to the largest residuals are selected for augmentation.

\textbf{Initial sampling.} To initialize the inner iteration of the low-rank SI, an initial set of $p$ angular directions $\Omega_{s^{(0)}_1}, \dots, \Omega_{s^{(0)}_p}$ is selected at random. Problem-specific initialization strategies may further enhance efficiency.

\begin{rem}
As an alternative to our proposed strategy, one may consider a sampling strategy similar to \cite{dektor2021rank} based on the discrete empirical interpolation method (DEIM) \cite{sorensen2016deim}.
Specifically, reduced solutions for \emph{all} unsampled angles can be computed via Galerkin projections in \eqref{eq:galerkin-projection}, yielding the low-rank approximation
\begin{equation}
\label{eq:lr_for_allangles}
\BPsi \approx \BX^{(n,k)}
\begin{bmatrix}
\bc^{(n,k)}_1 & \dots & \bc^{(n,k)}_{N_{\Omega}}
\end{bmatrix}
\triangleq
\BX^{(n,k)} \BC^{(n,k)}.
\end{equation}
For previously sampled angles, $\bc^{(n,k)}_j$ coincides with the corresponding column of $\BR^{(n,k)}$ in \eqref{eq:qr}. The matrix $\BC^{(n,k)} \in \mathbb{R}^{r_k \times N_{\BOmega}}$ collects the reduced coefficients of the angular flux in the current spatial basis. An angular basis and a coefficient core can then be obtained via SVD
\[
\BC^{(n,k)} = \BU^{(n,k)} \BS^{(n,k)} (\BV^{(n,k)})^T.
\]
New angles are subsequently selected by applying DEIM-induced sampling \cite{sorensen2016deim} to the reduced angular basis $\BV^{(n,k)}$.

In challenging benchmark tests, however, we observe that this DEIM-based strategy does not outperform the proposed subsampling approach in terms of sampling quality or the total number of augmented directions. Moreover, it incurs substantially higher computational cost, since Galerkin projections must be applied for all unsampled angles, rather than only for a small candidate set as in the proposed strategy.
\end{rem}

\begin{rem}
\review{
Our sampling strategy can be viewed as a combination of randomized and greedy sampling, similar in spirit to the randomized greedy algorithms developed in the reduced-order modeling community \cite{hesthaven2014efficient}. The randomized selection of candidate directions is introduced to improve efficiency, while the subsequent residual-based greedy selection within the candidate set exploits PDE information to identify the most informative samples and thereby improve sampling quality.
}

\review{
We also investigated a fully greedy strategy that selects each new sample by evaluating the residual over all unsampled angles. Though this approach can slightly reduce the total number of sampled directions compared to the proposed randomized greedy strategy, it substantially increases the sampling cost due to the significantly larger number of reduced solves and residual evaluations. In many of our numerical experiments, this overhead may make low-rank solver slower than the full-rank counterpart.}

\review{We also note that greedy angular sampling has been used to construct data-driven reduced-order models for RTE \cite{tencer2017accelerated,peng2022reduced}, whereas the present work focuses on designing offline-free rank-adaptive low-rank solvers.}

\end{rem}

\subsection{Stopping criteria for the inner-loop of low-rank SI\label{sec:stop}}
A commonly used stopping criterion in low-rank iterative solvers (see, e.g., \cite{appelo2025lraa}) measures the difference between two successive low-rank iterates, namely $$\|\BX^{(n,k)}\BS^{(n,k)}(\BV^{(n,k)})^T-\BX^{(n,k-1)}\BS^{(n,k-1)}(\BV^{(n,k-1)})^T\|_F\leq \textrm{threshold},$$
which can be computed through low-rank truncated summation without explicitly reconstructing the full matrices. Nevertheless, the procedure \cite{bachmayr2023low,einkemmer2025review} enlarges the low-rank space to dimension $r_k+r_{k-1}$ and typically requires a QR or SVD. For the challenging benchmarks considered in Sec. \ref{sec:motivation}, such augmentation can lead to substantial and unnecessary computational overhead.

To design an effective and efficient stopping criterion, we propose the following strategy, which avoids performing expensive truncated summation at every iteration.

\begin{enumerate}
\item When the low-rank SI starts, initialize 
\[
\bphi^{(n)}_{\mathrm{old}} := \bphi^{(n-1)}.
\]

\item In the $k$-th inner iteration of the low-rank SI, we first check whether the residuals for all $q$ candidate angles in the greedy random subsampling step are sufficiently small:
\begin{equation}
\max_{\Omega_j \in \mathcal{S}^{(k)}} 
\|\br_j^{(k)}\|_2
< \epsilon_{\mathrm{res}}.
\label{eq:residual-chk}
\end{equation}

\item If \eqref{eq:residual-chk} holds, we solve the Galerkin projection problem \eqref{eq:galerkin-projection} for all previously unsampled angles, yielding the low-rank approximation of the angular flux as in \eqref{eq:lr_for_allangles}.
The scalar flux can then be evaluated at cost 
$O(r_k (N_{\bx} + N_{\Omega}))$ as
\begin{equation}
\bphi^{(n,k)}
=
\BX^{(n,k)} \bigl( \BC^{(n,k)} \bw \bigr),
\qquad
\bw = (\omega_1,\dots,\omega_{N_{\Omega}})^T,
\end{equation}
where $\bw$ collects the normalized quadrature weights satisfying 
$\sum_{j=1}^{N_{\Omega}} \omega_j = 1$.

We then check whether
\begin{equation}
\|\bphi^{(n,k)} - \bphi^{(n)}_{\mathrm{old}}\|_2
\le
\epsilon_{\mathrm{diff}}.
\label{eq:phi-diff}
\end{equation}

If \eqref{eq:phi-diff} is not satisfied, update
\[
\bphi^{(n)}_{\mathrm{old}} := \bphi^{(n,k)}.
\]
Otherwise, terminate the inner loop and set
\[
\bphi^{(n,*)} = \bphi^{(n,k)}.
\]
\end{enumerate}

\textbf{Rationale for the stopping criterion.} The proposed stopping strategy combines two criteria: the residuals of the SI equation for candidate angles in the sampling step and the change in the scalar flux across inner iterations.
Note that relying solely on the residuals of candidate angles may lead to premature termination of the inner-loop iterations, which can adversely affect the convergence of the outer SI–DSA iterations. On the other hand, evaluating residuals for all angular directions incurs a prohibitive computational cost of \(O(N_{\Omega} N_{\bx} r_k)\). Meanwhile, using only the difference in scalar flux between successive iterations is also insufficiently robust. We numerically observe that, for certain angular directions, the residuals of the SI equation can remain significantly larger than the corresponding change in the scalar flux. In such cases, stopping based solely on the scalar-flux difference may cause premature termination of the inner loop as well, leading to stagnation in the outer-loop convergence.
In contrast, the proposed combined criteria achieve both efficiency and robustness across all benchmark tests. The resulting low-rank SI–DSA converges with essentially the same number of outer iterations as its full-rank counterpart.

\textbf{Construction of $\BV^{(n)}$ and $\BS^{(n)}$ after convergence.}
After the inner-loop iterations of the low-rank SI have converged, the angular basis and coefficient core can be obtained via a truncated SVD
$$
\BC^{(n,k)}
\approx
\BU^{(n)} \BS^{(n)} (\BV^{(n)})^T,
$$
where
$\BU^{(n)} \in \mathbb{R}^{r_k \times r'}$,
$\BS^{(n)} \in \mathbb{R}^{r' \times r'}$,
$\BV^{(n)} \in \mathbb{R}^{N_{\Omega} \times r'}$,
and $r' \le r_k$.
The spatial basis is then updated as
\[
\BX^{(n)} = \BX^{(n,k)} \BU^{(n)}.
\]
This truncation step is optional, since only the scalar flux $\bphi^{(n,*)}$ is required to proceed with the outer SI–DSA iterations. Even when the angular flux is needed, for instance, in implicit time marching of time-dependent problems, a single truncation performed after  convergence of the outer loop is sufficient.

\begin{rem}
\review{
We emphasize that, upon convergence, Galerkin projections are applied only to the unsampled angular directions to approximate their reduced coordinates. For the sampled directions, the corresponding SI equations \eqref{eq:SI} have already been solved directly by transport sweeps, and their reduced coordinates have been computed during the incremental basis update.}
\end{rem}

\review{
\begin{rem}
(Guidelines for choosing thresholds.)
The tolerances are chosen based on two considerations. First, the inner-loop thresholds should be smaller than the outer SI-SA tolerance $\epsilon_{\textrm{SI-SA}}$ to prevent premature termination of the inner iterations, which may slow down or even hinder the outer-loop convergence. Second, the truncation tolerance in SVD, namely $\epsilon_{\textrm{SVD}}$, should be consistent with the inner stopping thresholds, since it is introduced primarily as a stabilization mechanism to eliminate components whose magnitudes fall below the prescribed accuracy.  In our tests, we use
\begin{equation}
\epsilon_{\textrm{SVD}} \lesssim \epsilon_{\textrm{diff.}},\epsilon_{\textrm{res.}}
\quad\text{and}\quad
\epsilon_{\textrm{diff.}},\epsilon_{\textrm{res.}} \approx 0.1\epsilon_{\textrm{SI-SA}}.
\end{equation}
Adaptive threshold control, similar to \cite{bachmayr2017iterative}, may further improve efficiency but defers to future investigation.
\end{rem}
}
\subsection{Incremental basis and reduced operator update\label{sec:lr-update}}
In the inner-loop iteration of our low-rank SI, the space of spatial modes is progressively enlarged by adding snapshots of angular fluxes for sampled angles. By exploiting this hierarchical structure, the spatial basis $\BX^{(n,k)}$ and Galerkin projected operators onto its span can be updated incrementally. Specifically, the orthonormal spatial basis $\BX^{(n,k)}$ is updated via \review{a variant of modified Gram--Schmidt (MGS) procedure with adaptive reorthogonalization  \cite{daniel1976reorthogonalization}, as summarized in Alg. \ref{alg:mgs2-a}. Since we apply at least two MGS projections, we call it MGS2-A.}

\review{It is well-known that orthonormal bases obtained by MGS may gradually lose orthogonality due to roundoff accumulation, thereby reducing the accuracy and numerical stability of the projected operators. Instead of applying an expensive full QR to the accumulated basis, MGS2-A effectively mitigates this loss of orthogonality through additional projections that preserve the hierarchical basis structure \cite{daniel1976reorthogonalization}.} 

\review{When $p>1$, the sampled snapshots are processed sequentially: after each sampled snapshot is appended, it immediately becomes part of the active basis used to orthogonalize the remaining snapshots in the same batch.}

\review{The MGS2-A algorithm preserves the hierarchical structure of basis, namely,}
${\small{\BX^{(n,k)}=\begin{bmatrix} \BX^{(n,k-1)} & \BP^{(n,k)} \end{bmatrix}}}$ with $\BP^{(n,k)}\in\mathbb{R}^{N_{\bx}\times p}$. Then, the projected discrete total cross section can be incrementally updated as  (dropping the superscript index 
$(n,k-1)$ for brevity):
\begin{align}
(\BX^{(n,k)})^T\BSigma_t\BX^{(n,k)}=(\begin{bmatrix} \BX & \BP \end{bmatrix})^T\BSigma_t \begin{bmatrix} \BX & \BP \end{bmatrix}=\begin{bmatrix}
\underbrace{\BX^T\BSigma_t\BX}_{\textrm{Known}} & \BX^T(\BSigma_t\BP)\\
(\BP^T\BSigma_t)\BX & \BP^T\BSigma_t\BP
\end{bmatrix}.
\end{align}
The projections of other operators can be updated in the same manner. In particular, the projection of the transport operator for the
angle $\BOmega_j$, given by
$(\BX^{(n,k)})^T \BD_j \BX^{(n,k)}$, is  a linear combination of
the projections of $\BD_x^{\pm}$ and $\BD_y^{\pm}$ defined in
\eqref{eq:upwind-operator}.

\begin{algorithm}[ht]
\caption{\review{Double modified Gram--Schmidt projections with adaptive reorthogonalization  (MGS2-A)}\label{alg:mgs2-a}}
\begin{algorithmic}[1]
\STATE{\textbf{Input}: Orthonormal spatial basis $\BX^{(n,k-1)}$, coefficient matrix $\BR^{(n,k-1)}$, angular snapshots $\bpsi^{(n)}_{\sk_1},\dots,\bpsi^{(n)}_{\sk_p}$, tolerance $\epsilon_{\mathrm{MGS}}$, and maximum cleanup count $m_{\max}$.}
\STATE{Set $\BX:=\BX^{(n,k-1)}$ and $\BR:=\BR^{(n,k-1)}$.}
\FOR{$l=1:p$}
    \STATE{\review{Set $\hat{\bpsi}:=\bpsi^{(n)}_{\sk_l}$ and $\bc:=\mathbf{0}$.}}
    \FOR{\review{$t=1:2$}}
        \STATE{\review{$\Delta\bc=\BX^T\hat{\bpsi}$,}\quad\review{$\hat{\bpsi}\leftarrow \hat{\bpsi}-\BX\Delta\bc$,}\quad\review{$\bc\leftarrow \bc+\Delta\bc$.}}
    \ENDFOR
    \STATE{\review{$\alpha=\|\hat{\bpsi}\|_2$ and $\hat{\bpsi}\leftarrow \hat{\bpsi}/\alpha$.}}
    \FOR{\review{$m=1:m_{\max}$}}
        \IF{\review{$|\BX(:,1)^T\hat{\bpsi}|\le \epsilon_{\mathrm{MGS}}$}}
            \STATE{\review{\textbf{break}}}
        \ENDIF
        \STATE{\review{$\hat{\bpsi}\leftarrow \alpha\hat{\bpsi}$,}\quad \review{$\Delta\bc=\BX^T\hat{\bpsi}$.}}
        \STATE{\review{$\hat{\bpsi}\leftarrow \hat{\bpsi}-\BX\Delta\bc$,}\quad\review{$\bc\leftarrow \bc+\Delta\bc$.}}
        \STATE{\review{$\alpha=\|\hat{\bpsi}\|_2$ and $\hat{\bpsi}\leftarrow \hat{\bpsi}/\alpha$.}}
    \ENDFOR
    \IF{\review{$|\BX(:,1)^T\hat{\bpsi}|>\epsilon_{\mathrm{MGS}}$}}
        \STATE{\review{Record one MGS orthogonality failure.}}
    \ENDIF
    \STATE{\review{Update $\BX:=\begin{bmatrix}\BX & \hat{\bpsi}\end{bmatrix}$, define $\tilde{\bc}:=\bigl(\bc^T,\alpha\bigr)^T$, and set $\BR:=\begin{bmatrix}\BR & \tilde{\bc}\end{bmatrix}$.}}
\ENDFOR
\STATE{\review{Set $\BX^{(n,k)}:=\BX$ and $\BR^{(n,k)}:=\BR$.}}
\end{algorithmic}
\end{algorithm}

\subsection{Computational cost of each inner-loop iteration in low-rank SI}
Here, we summarize the computational cost of the $k$-th inner-loop iteration in low-rank SI. 
\begin{enumerate}
    \item \textbf{Transport sweeps}: Solving the SI equation \eqref{eq:SI} for the $p$ newly sampled directions requires $O(p N_{\mathbf{x}})$ operations.
    \item \textbf{Basis and operators update}: Without \review{full QR}, the incremental update of the spatial basis and projected operators incurs a cost of $O(N_{\mathbf{x}} r_{k-1} p)$. 
    \item \textbf{Greedy random subsampling}: Solving the projected systems for the $q$ candidate angles and computing their corresponding residuals incur a cost of $O(q r_k^3)$ and $O(N_{\mathbf{x}} r_k q)$, respectively.
    \item \textbf{Convergence monitoring}: The cost of evaluating stopping criteria depends on the residual magnitude. If the maximum residual among candidate angles exceeds the tolerance, no further cost is incurred. Otherwise, generating the coefficient matrix $\mathbf{C}^{(n,k)}$ requires additional $O((N_{\mathbf{\Omega}} - q - r_k) r_k^3)$ operations, and updating the scalar flux costs $O((N_{\mathbf{x}} + N_{\BOmega}) r_k)$.
\end{enumerate}
We emphasize that, as observed in our numerical experiments,  projections for the full set of unsampled angular directions, together with the subsequent angular integration, are performed only in the final few iterations before convergence, thereby significantly reducing the average per-iteration cost. 

\subsection{Contributions and relation to prior work\label{sec:contribution}}

Fixed-rank sweeping-based DLR methods for the first-order formulation of the RTE have been developed in \cite{peng2023sweep,haut2026efficient}, \review{and could potentially be integrated with rank-adaptive BUG \cite{ceruti2022rank}.} \review{The main contribution of this work is not merely the introduction of rank adaptivity, but the development of a novel and efficient rank-adaptation strategy that avoids $O(r)$ space augmentation required by classical  approaches.} 

Although both our method and the concurrent $S_N$-like collocation-based DLR method in \cite{haut2026efficient} construct low-rank approximations via angular sampling and physical-space projections, the underlying algorithmic structures differ fundamentally. The method in \cite{haut2026efficient} is designed for implicit time marching of time-dependent thermal radiation transfer problems. At each time step, \review{angles are sampled} using a DEIM-based sampling procedure applied to the low-rank angular basis from the previous time step, followed by the application of SI–DSA to a modified $S_N$ system determined by those angles. Thus, angular sampling is performed outside the SI–DSA iteration in \cite{haut2026efficient}. In contrast, our method integrates adaptive angular sampling directly into each inner SI iteration. \review{Unlike time-dependent settings, steady-state problems lack previous-time-step solutions to guide sampling, while the intermediate solutions may vary substantially throughout the iterative solve. Consequently, sampling should be performed adaptively within the SI iteration. Our randomized, residual-informed sampling strategy provides this capability and is the key mechanism enabling rank adaptivity with only mild space augmentation.}

\review{From an implementation perspective, a key advantage of our formulation is that it preserves the standard one-angle transport sweep. Let $s$ denote the number of spatial DOFs per element. For each sampled angular direction, our method retains the original $s\times s$ local system arising from the DG discretization, allowing the reuse of existing optimized one-angle transport sweep implementations without modification. The additional low-rank operations enter through angle sampling, physical-space projections and basis augmentation. This advantage is also shared by the $S_N$-like DLR method in \cite{haut2026efficient}.  In contrast, the method in \cite{peng2023sweep} performs projections in the angular space, requiring the transport sweep to solve a modified $rs\times rs$ local system in each element and thereby requiring intrusive modifications to transport sweep implementations.}

\review{In summary, the DEIM-based approach in \cite{haut2026efficient} is designed for implicit time marching and exploits previous-time-step angular basis to sample angles outside SI--DSA, enabling the direct reuse of existing SI--DSA code straightforward. In contrast, our approach is tailored to steady-state problems and achieves efficient rank adaptivity by integrating angular sampling and low-rank operations into  SI. These two methods are complementary and can be naturally  combined in implicit time marching: the DEIM-based approach can provide high-quality initial angular samples, while our method can efficiently adapt the solution rank during SI if necessary. This combination may be attractive when large time steps reduce the predictive quality of the previous-time-step basis.}

To further place our work in the broader context of efficient rank adaptation, we note that aggressive  augmentation was also mitigated in \cite{hu2022adaptive} and later integrated into spatial domain decomposition in \cite{brunner2026domain}. In that approach, the basis is enriched at the beginning of each time step with random orthogonal vectors whose number equals to the rank of the boundary conditions. In contrast, our method targets iterative solvers and enriches the basis through error-indicated sampling rather than  random padding.

\subsection{Future extensions to  implicit time marching\label{sec:extension}}
The proposed low-rank SI-DSA can be straightforwardly extended to solve linear systems arising from implicit time marching, as they can be viewed as modified steady-state problems in each time step. For example, each time step of the backward Euler scheme can be interpreted as solving a steady-state problem with a modified discrete total cross section $\widetilde{\BSigma}_t=\BSigma_t+\frac{1}{\Delta t}\BM$, where $\BM$ denotes a mass matrix and $\Delta t$ is the time step size. 
As the primary focus of this work is the development  of efficient iterative solvers, a detailed investigation of such time-dependent extensions is deferred to future work.

\section{Numerical examples\label{sec:numerical}}
We demonstrate the performance of our method via a series of benchmark tests \review{in the 2D X-Y geometry}. Throughout this section ``LR" and ``FR" refer to low-rank and full-rank solutions, respectively.

In all tests, we use piecewise linear polynomial basis for the DG discretization. A fully-consistent DSA preconditioner is employed \cite{Adams2002FastIM} (see Appx. B of \cite{peng2024reduced} for details). \review{The code was implemented in {\tt Matlab}.} The diffusion correction equation in DSA is solved using the conjugate gradient (CG) method preconditioned by a V-cycle algebraic multigrid (AMG) method with $10^{-12}$ as the relative residual tolerance. \review{This CG-AMG solver is implemented using the {\tt iFEM} package \cite{chen2009integrated}}. For comparison, \review{unless otherwise specified,} the stopping criteria is chosen as $\|\bphi^{(n)}-\bphi^{(n-1)}\|_\infty\leq\epsilon_{\textrm{SI-SA}}=10^{-6}$ for both the full-rank SI-DSA and the outer-loop iteration of low-rank SI-DSA (see Alg. \ref{alg:lr-si-dsa} for details). In addition, a zero initial guess is used for both methods. The initial angle samples for the first  inner-loop iteration of low-rank SI-DSA are chosen randomly.

In our low-rank SI-DSA, unless otherwise specified, we set $\epsilon_{\textrm{res.}}=\epsilon_{\textrm{diff.}}=10^{-7}$ in the stopping criteria for the inner-loop iterations of low-rank SI, and choose $p=1,q=8$ in the residual-based greedy random subsampling of angles. In other words, by default, we select one angle sample from $8$ random candidates. 
In addition, we construct the low-rank core and the angular flux by truncating at the rank $r$ satisfies $r=\arg\min_l \sum_{i=1}^l s_i/\sum_{i=1}^{r_n} s_i\leq \epsilon_{\textrm{SVD}}$, where $s_i$ are singular values in the SVD of $C^{(n)}$ after the inner-loop convergence and  $\epsilon_{\textrm{SVD}}=10^{-8}$. \review{When updating spatial basis, we apply at least $2$ and at most $6$ MGS projections and set $\epsilon_{\textrm{MGS}}=10^{-10}$ in Alg. \ref{alg:mgs2-a}. In all tests, the orthogonality test in Alg. \ref{alg:mgs2-a} of spatial basis eventually always passes.}

\review{All the experiments were performed using {\tt Matlab R2024a} with {\tt OpenBLAS 0.3.24} and {\tt LAPACK 3.9.1} backends on a 2021 24-inch {\tt iMac} with a $16$ GB memory {\tt Apple M1 chip}. The wall-clock times  were measured using the {\tt tic-toc} function.}

To quantify the performance of our method, we define metrics including the speedup, compression ratio and inner-loop oversampling ratio as follows.
\begin{align}
\textrm{speedup}&=\frac{\textrm{Computational Time of \review{Full-rank} SI-DSA}}{\textrm{Computational Time of \review{Low-rank} SI-DSA}},\\
\textrm{compression ratio}&=\frac{\textrm{DOFs for $\bpsi$ in Low-rank SI-DSA}}{\textrm{DOFs for $\bpsi$ in Full-rank SI-DSA}},\\
\textrm{inner-loop oversampling ratio}&=\frac{\textrm{Number of Sampled Angles}-\textrm{Solution Rank}}{\textrm{Solution Rank}}.
\end{align}

Before delving into details, we outline numerical tests performed and their purposes. In Sec. \ref{sec:homo}, for a homogeneous medium ranging from weak to strong scattering, we perform a refinement study to investigate the influence of mesh resolution and the effectiveness of our method in different regimes. 
\pzc{In Sec. \ref{sec:lattice}, a lattice problem is solved  to illustrate the 
\review{benefits of using MGS2-A algorithm in Alg. \ref{alg:mgs2-a}} to incrementally update spatial basis.}
In Sec. \ref{sec:variable-scattering}, we simulate a multiscale variable scattering problem and investigate the influence of $p$ and $q$ in our sampling step. In Sec. \ref{sec:pin-cell}, a pin-cell problem with strong material discontinuities is solved repeatedly to demonstrate the robustness of the method with respect to randomness in angular sampling. 

\subsection{Homogeneous medium in different regime\label{sec:homo}}
\begin{figure}[]
\centering
    \includegraphics[width=0.45\textwidth]{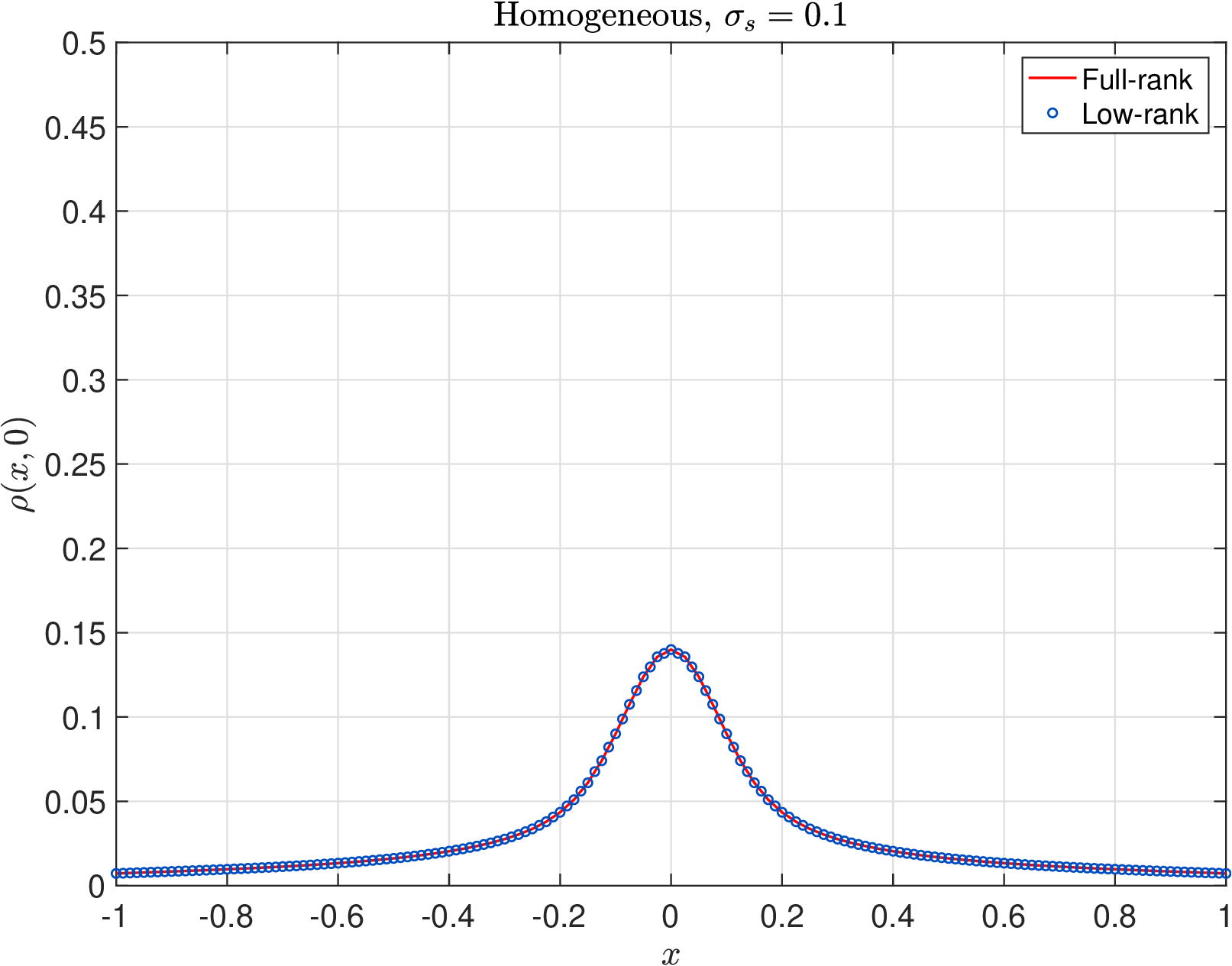}
    \includegraphics[width=0.45\textwidth]{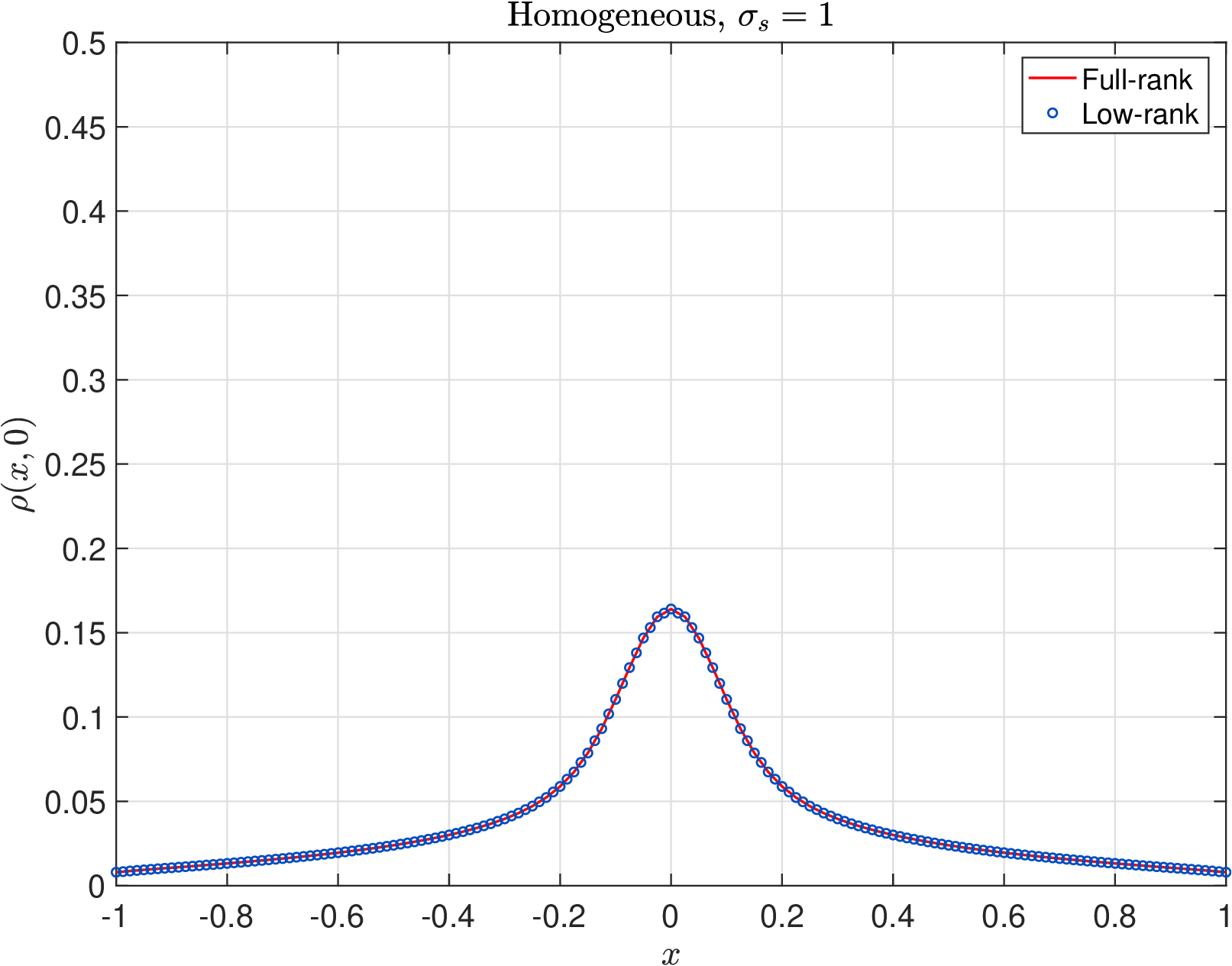}
    \includegraphics[width=0.45\textwidth]{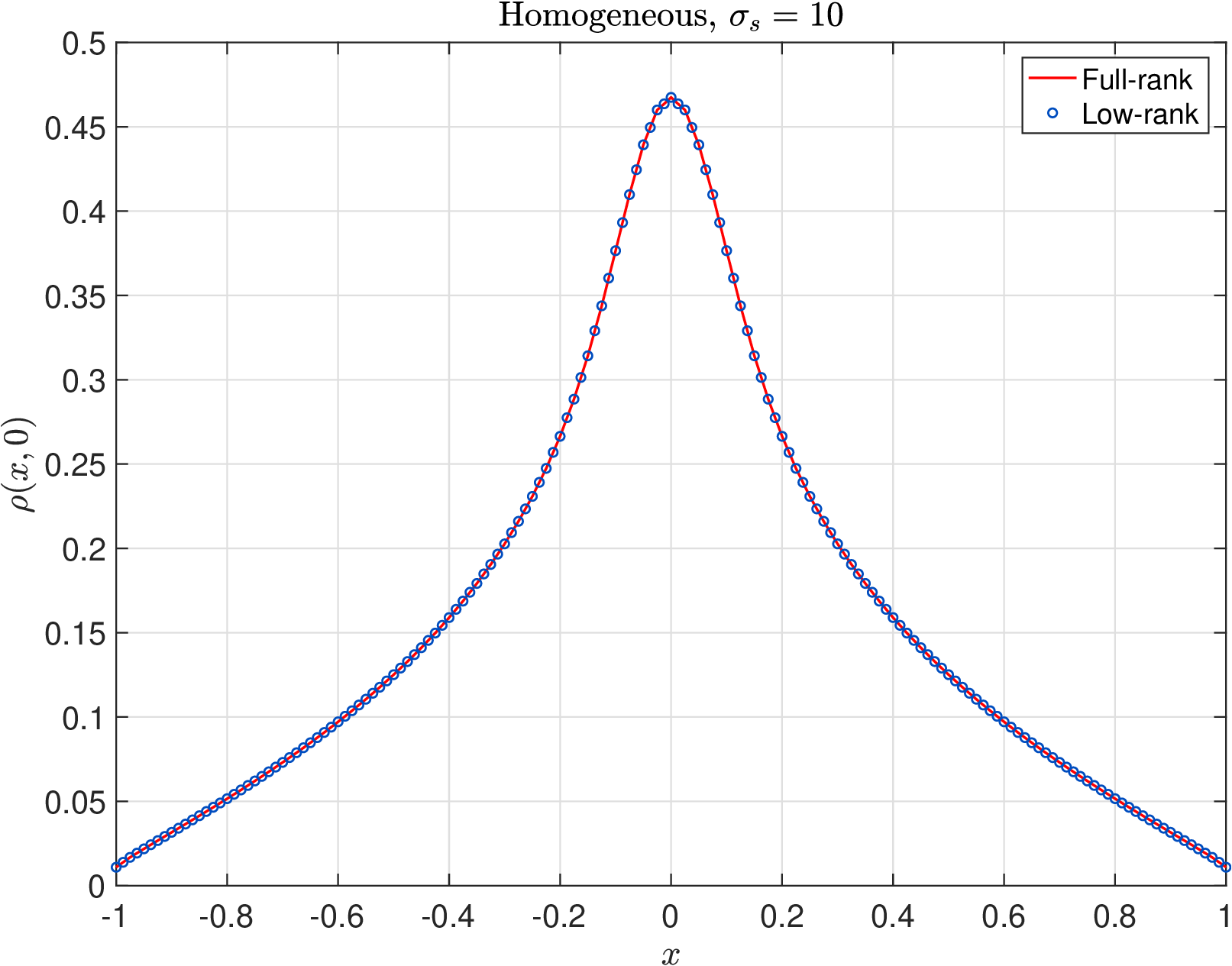}
    \includegraphics[width=0.45\textwidth]{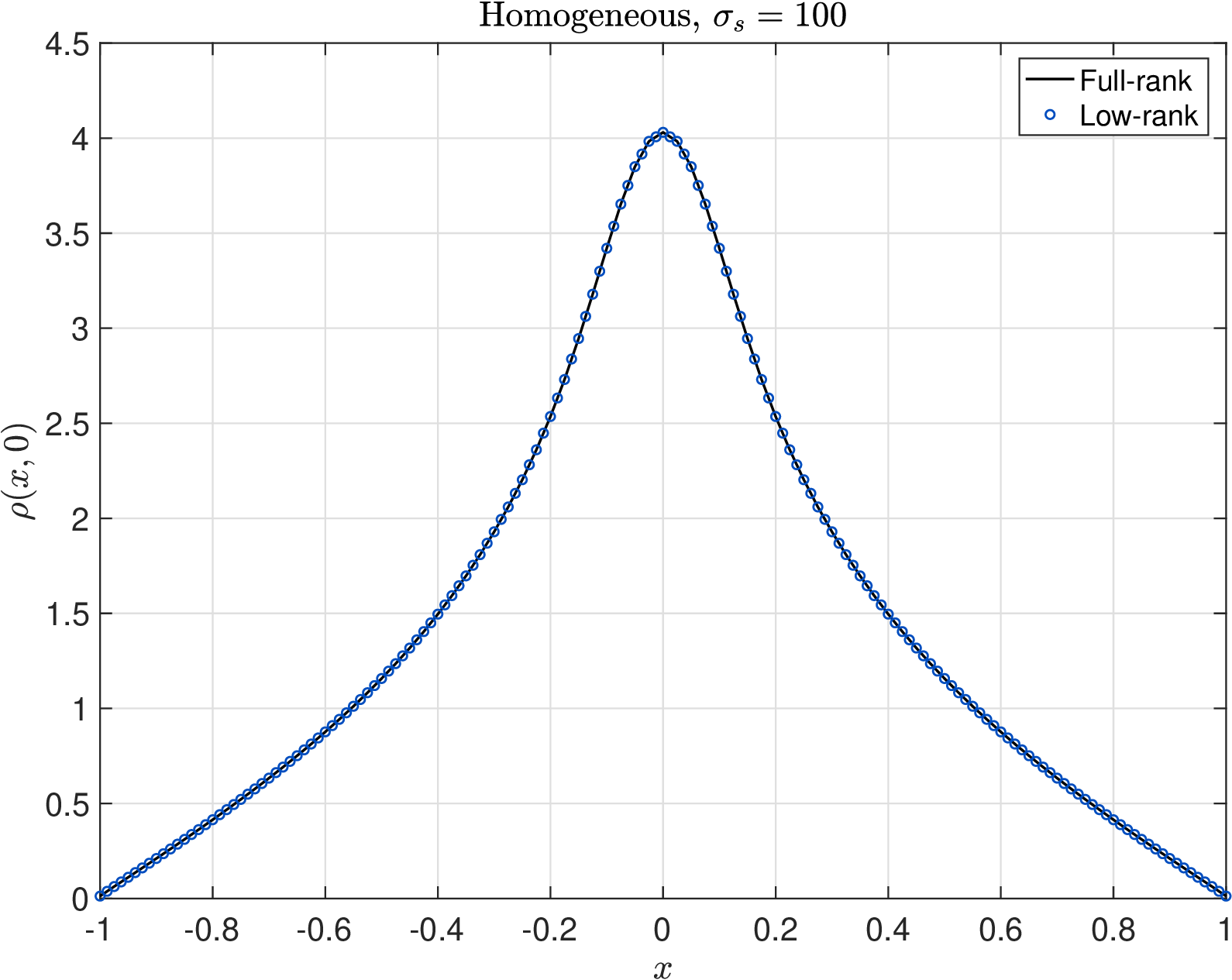}
  
  \caption{\review{Comparison between the scalar flux for the full-rank and low-ran solutions} of the homogeneous problem  in Sec. \ref{sec:homo} \review{along $y=0$} with $(N_x,N_y,N_\theta,N_{\Omega_z})=(80,80,40,20)$. Top left: $\sigma_s=0.1$. Top right: $\sigma_s=1$. Bottom left: $\sigma_s=10$. Bottom right: $\sigma_s=100$. \label{fig:homo}}
\end{figure}

We consider the computational domain $[-1,1]^2$ with no absorption, $\sigma_s(x,y)=0.1,1,10,100$ and zero inflow boundary conditions. An isotropic Gaussian source $G(x,y)=\exp(-100(x^2+y^2))$ is imposed.  We perform a refinement study using a $16L\times 16L$ rectangular mesh in the physical space and the CL($8L,4L$) quadrature in the angular space for $L=2,3,4,5$. \review{We compare the 1D profile of the full-rank and low-rank solutions along $y=0$ with $L=5$} in Fig. \ref{fig:homo}, and they visually match each other well.

The main goal of this example is to investigate the influence of $\sigma_s$ and the mesh resolution on our low-rank SI-DSA. The
results are summarized in Tab. \ref{tab:homo}. Main observations are as follows. 

\textbf{Accuracy and number of iterations for convergence.} Regardless of the mesh resolution and the strength of the scattering effect, both the low-rank SI-DSA and full-rank SI-DSA converge within exactly the same number of iterations. Moreover, the difference between the scalar fluxes computed by low-rank and full-rank SI-DSA is at most on the order of $10^{-7}$, which is at least one order of magnitude smaller than $\epsilon_{\textrm{SI-SA}}$ in the outer-loop stopping criterion of SI-DSA. These observations demonstrate the correctness and effectiveness of our low-rank SI-DSA.

\textbf{Influence of mesh resolution.} For all $\sigma_s$, increasing the mesh resolution leads to greater memory savings, and the speedup of the low-rank SI-DSA over the full-rank counterpart generally increases as well.

\textbf{Influence of scattering strength $\sigma_s$.} 
\begin{enumerate}
    \item \textbf{Convergence.} When $\sigma_s=0.1$ (transport regime) or $100$ (diffusion regime), both full-rank and low-rank SI-DSA converges within $4-5$ iterations. In contrast, more iterations are required for convergence in the intermediate regime ($\sigma_s=1,10$).
    \item \textbf{Rank \review{and compression ratio}.} The solution rank almost grows linearly with respect to the grid resolution in the transport regime ($\sigma_s=0.1$), scales roughly as $O(L^{1.5})$ in the intermediate regime ($\sigma_s=1,10$), and almost remains constant once the resolution is sufficiently fine in the diffusion regime ($\sigma_s=100$).
    \review{Since the DOFs in the physical and angular space scales as $L^2$, the effective DOFs scales as  $L^2r$ in the low-rank solvers with $r$ being the rank. As a result, the compression ratio scales as $L^{-1}$ in the transport regime, $L^{-0.5}$ in the intermediate regime and $L^{-2}$ in the diffusion regime.} The resolution dependence of the rank in the transport and intermediate regimes is consistent with observations in \review{Sec. 4.1 and Sec. 4.2} of our previous work \cite{guo2025inexact}, \review{and more detailed tests and explanations can be found there}.
\end{enumerate}

    \textbf{Computational savings compared to full-rank solver.} The low-rank SI-DSA consistently achieves computational time and memory savings over the full-rank SI-DSA across all configurations in this test. The largest savings occur in the diffusion regime  ($\sigma_s=100$). In particular, when $\sigma_s=100$ and $L=5$, our low-rank SI-DSA gains more than \review{$10\times$} reductions in both the computational time and memory.  
\begin{table}[htbp]
  \centering
 \medskip
  \begin{subtable}{\textwidth}
    \small
   \centering
    \begin{tabular}
    {|l|c|c|c|c|c|c|c|}
    \hline
    $L$ &  Rank & $\|\bpsi_{\textrm{LR}}-\bpsi_{\textrm{FR}}\|_2$ &$\|\bphi_{\textrm{LR}}-\bphi_{\textrm{FR}}\|_2$ &  C-R & FR-Iter & LR-Iter & Speedup    \\ \hline
$ \review{2} $ & $ \review{53} $ & $\review{5.53 \times 10^{-8}}$ & $\review{1.16 \times 10^{-8}}$ & $\review{42.71\%}$ & $\review{5}$ & $\review{5}$ & $\review{1.72\times}$ \\ 
$ \review{3} $ & $ \review{95} $ & $\review{6.03 \times 10^{-9}}$ & $\review{1.91 \times 10^{-10}}$ & $\review{34.02\%}$ & $\review{5}$ & $\review{5}$ & $\review{2.76\times}$ \\ 
$ \review{4} $ & $ \review{128} $ & $\review{1.06 \times 10^{-7}}$ & $\review{3.65 \times 10^{-9}}$ & $\review{25.78\%}$ & $\review{5}$ & $\review{5}$ & $\review{2.93\times}$ \\ 
$ \review{5} $ & $ \review{164} $ & $\review{5.67 \times 10^{-7}}$ & $\review{1.41 \times 10^{-8}}$ & $\review{21.14\%}$ & $\review{5}$ & $\review{5}$ & $\review{3.40\times}$ \\ \hline
    \end{tabular}
      \caption{$\sigma_s=0.1$.}
\end{subtable}
\\
\begin{subtable}{\textwidth}
    \small
   \centering
    \begin{tabular}
    {|l|c|c|c|c|c|c|c|}
    \hline
    $L$ &  Rank & $\|\bpsi_{\textrm{LR}}-\bpsi_{\textrm{FR}}\|_2$ & $\|\bphi_{\textrm{LR}}-\bphi_{\textrm{FR}}\|_2$ & C-R & FR-Iter & LR-Iter & Speedup \\ \hline
$ \review{2} $ & $ \review{64} $ & $\review{2.17 \times 10^{-7}}$ & $\review{1.54 \times 10^{-7}}$ & $\review{51.57\%}$ & $\review{8}$ & $\review{9}$ & $\review{1.26\times}$ \\ 
$ \review{3} $ & $ \review{122} $ & $\review{1.57 \times 10^{-7}}$ & $\review{1.55 \times 10^{-8}}$ & $\review{43.69\%}$ & $\review{8}$ & $\review{8}$ & $\review{2.15\times}$ \\ 
$ \review{4} $ & $ \review{176} $ & $\review{4.52 \times 10^{-6}}$ & $\review{1.25 \times 10^{-7}}$ & $\review{35.45\%}$ & $\review{8}$ & $\review{8}$ & $\review{2.41\times}$ \\ 
$ \review{5} $ & $ \review{223} $ & $\review{2.75 \times 10^{-6}}$ & $\review{5.15 \times 10^{-8}}$ & $\review{28.75\%}$ & $\review{8}$ & $\review{8}$ & $\review{2.50\times}$ \\ \hline
    \end{tabular}
      \caption{$\sigma_s=1$}
\end{subtable}
\\
\begin{subtable}{\textwidth}
\small
   \centering
    \begin{tabular}
    {|l|c|c|c|c|c|c|c|}
    \hline
    $L$ &  Rank & $\|\bpsi_{\textrm{LR}}-\bpsi_{\textrm{FR}}\|_2$ & $\|\bphi_{\textrm{LR}}-\bphi_{\textrm{FR}}\|_2$ & C-R & FR-Iter & LR-Iter & Speedup \\ \hline
$ \review{2} $ & $ \review{64} $ & $\review{3.17 \times 10^{-7}}$ & $\review{2.91 \times 10^{-7}}$ & $\review{51.57\%}$ & $\review{8}$ & $\review{8}$ & $\review{1.39\times}$ \\ 
$ \review{3} $ & $ \review{139} $ & $\review{8.05 \times 10^{-8}}$ & $\review{2.49 \times 10^{-8}}$ & $\review{49.78\%}$ & $\review{8}$ & $\review{8}$ & $\review{1.83\times}$ \\ 
$ \review{4} $ & $ \review{191} $ & $\review{8.10 \times 10^{-8}}$ & $\review{1.10 \times 10^{-8}}$ & $\review{38.47\%}$ & $\review{8}$ & $\review{8}$ & $\review{2.31\times}$ \\ 
$ \review{5} $ & $ \review{212} $ & $\review{1.48 \times 10^{-7}}$ & $\review{8.46 \times 10^{-9}}$ & $\review{27.33\%}$ & $\review{8}$ & $\review{8}$ & $\review{3.16\times}$ \\ \hline
    \end{tabular}
      \caption{$\sigma_s=10$}
\end{subtable}
\\
\begin{subtable}{\textwidth}
\small
   \centering
    \begin{tabular}
    {|l|c|c|c|c|c|c|c|}
    \hline
    $L$ &  Rank & $\|\bpsi_{\textrm{LR}}-\bpsi_{\textrm{FR}}\|_2$ & $\|\bphi_{\textrm{LR}}-\bphi_{\textrm{FR}}\|_2$ & C-R & FR-Iter & LR-Iter & Speedup \\ \hline
$ \review{2} $ & $ \review{53} $ & $\review{7.17 \times 10^{-8}}$ & $\review{7.12 \times 10^{-8}}$ & $\review{42.71\%}$ & $\review{4}$ & $\review{4}$ & $\review{2.55\times}$ \\ 
$ \review{3} $ & $ \review{61} $ & $\review{4.51 \times 10^{-7}}$ & $\review{4.51 \times 10^{-7}}$ & $\review{21.84\%}$ & $\review{4}$ & $\review{4}$ & $\review{4.60\times}$ \\ 
$ \review{4} $ & $ \review{62} $ & $\review{2.15 \times 10^{-7}}$ & $\review{2.14 \times 10^{-7}}$ & $\review{12.49\%}$ & $\review{4}$ & $\review{4}$ & $\review{7.53\times}$ \\
$ \review{5} $ & $ \review{64} $ & $\review{1.39 \times 10^{-7}}$ & $\review{1.37 \times 10^{-7}}$ & $\review{8.25\%}$ & $\review{5}$ & $\review{5}$ & $\review{10.40\times}$ \\ \hline
    \end{tabular}
      \caption{$\sigma_s=100$}
\end{subtable}
\caption{Results for the homogeneous medium problem in Sec. \ref{sec:homo}. C-R: compression ratio. FR/LR-Iter: number of source iterations required by convergence for the full/low-rank SI-DSA. $L$: refinement level corresponding to   $(N_x,N_y,N_\theta,N_{\Omega_z})=(16L,16L,8L,4L)$.\label{tab:homo} }
\end{table}
\subsection{Lattice problem\label{sec:lattice}}
\begin{figure}[]
\centering
    \hspace{-0.5cm}
    \includegraphics[trim=10mm 9mm 10mm 0mm, clip,width=0.33\textwidth]{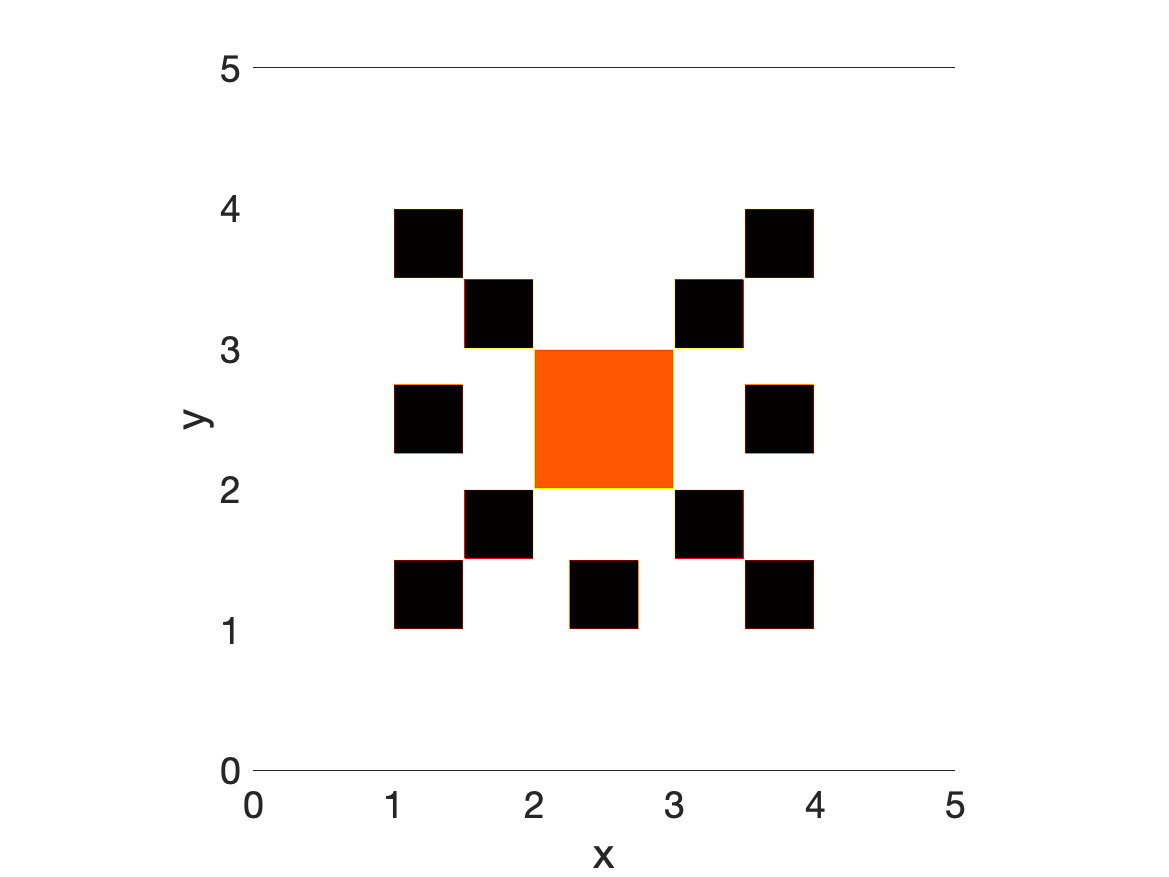}
    \hspace{-0.5cm}
    \includegraphics[trim=10mm 9mm 10mm 0mm, clip,width=0.33\textwidth]{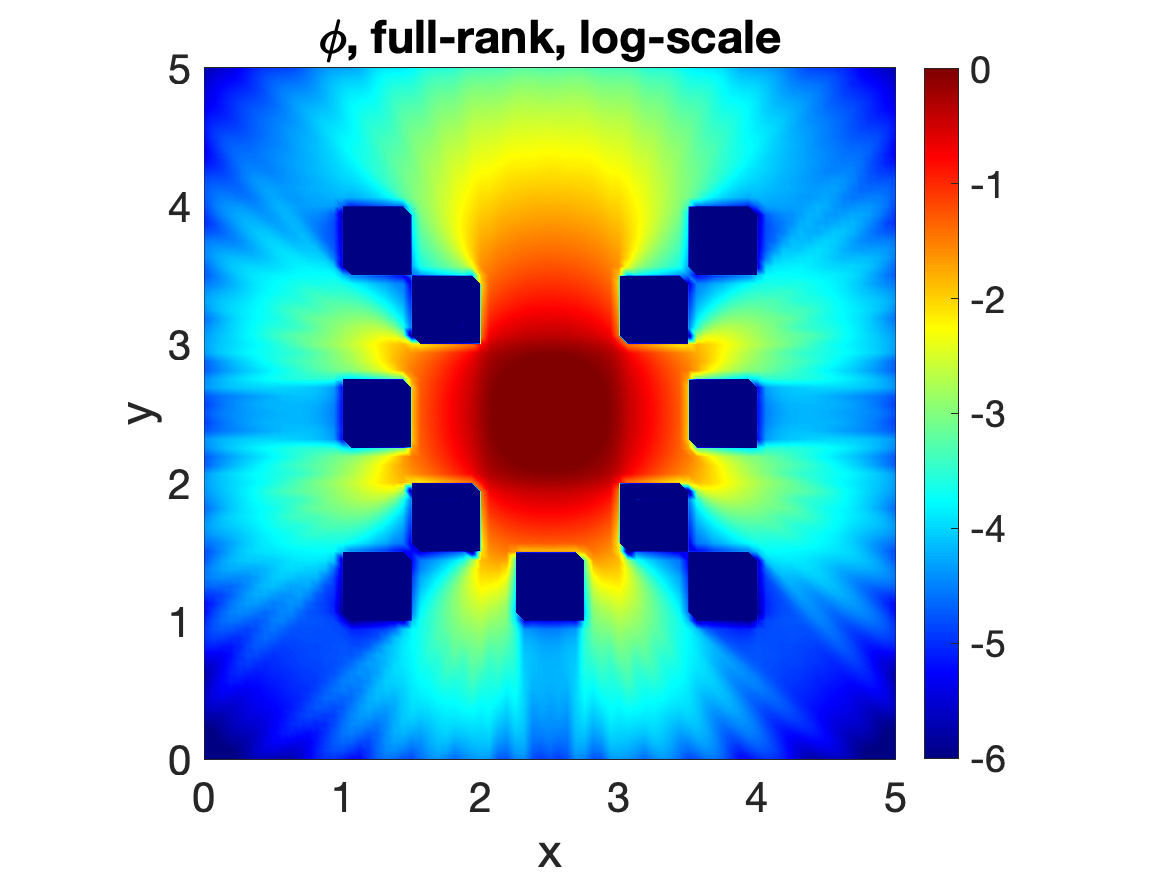}
    \hspace{-0.5cm}
    \includegraphics[trim=10mm 9mm 10mm 0mm, clip,width=0.33\textwidth]{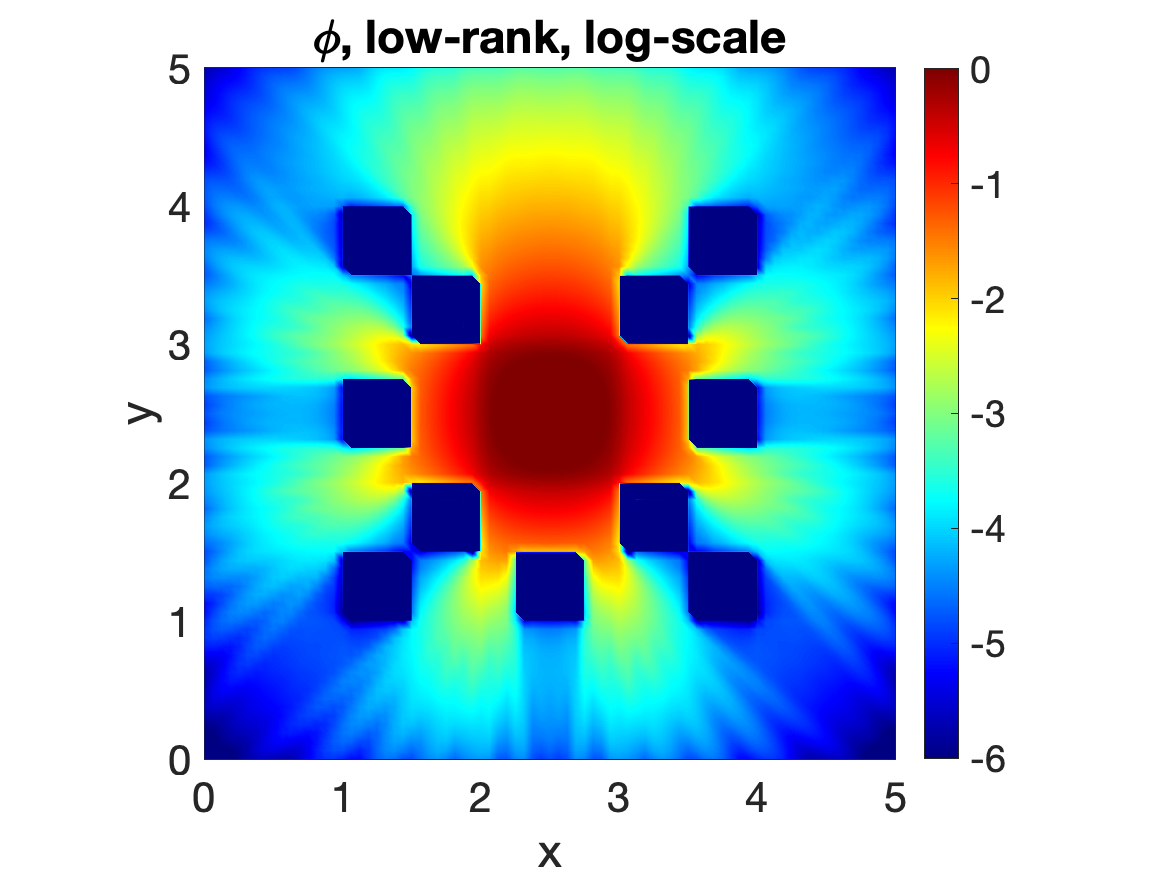}
  \caption{Configuration of scattering cross section and scalar fluxes obtained by full-rank and low-rank SI-DSA for the Lattice example in Sec. \ref{sec:variable-scattering}. We note that negative scalar flux is generated by upwind DG method in this example, and we round negative values to $10^{-16}$ for plotting purposes. Left: configuration of $\sigma_s(x,y)$, white and orange region pure scattering, black region pure absorption and a constant source imposed in the orange region.   Middle: $\phi$ obtained by full-rank SI-DSA. Right: $\phi$ obtained by low-rank SI-DSA. \label{fig:lattice}}
\end{figure}

\pzc{We solve} a lattice problem with zero inflow boundary conditions, and its configuration is shown in Fig. \ref{fig:lattice}. The black region corresponds to pure absorption with $(\sigma_a,\sigma_s)=(100,0)$, while the white and orange regions correspond to pure scattering with $(\sigma_a,\sigma_s)=(0,1)$. A constant source $G(x,y)=1$ is imposed in the orange region.
The computational domain $[0,5]^2$ is discretized using an $80\times80$ mesh in physical space and the CL(32,16) quadrature in angular space. The scalar fluxes obtained by the low-rank and full-rank SI-DSA are shown on a logarithmic scale in Fig. \ref{fig:lattice} and match each other well.

This example is used to demonstrate the necessity of \review{using MGS2-A in Alg. \ref{alg:mgs2-a}} to incrementally update spatial basis and the projected reduced operators. In particular, we compare the performance of our method with \review{direct full QR decomposition} and \review{single MGS projection using full QR for reorthogonalization when the orthogonality test fails.}  We note that the incremental update of projection operators relies on the hierarchical structure of the spatial basis; therefore, incremental updates of the projected operators cannot be applied \pzc{when using full QR}.

\begin{figure}[]
\centering
    \includegraphics[width=0.45\textwidth]{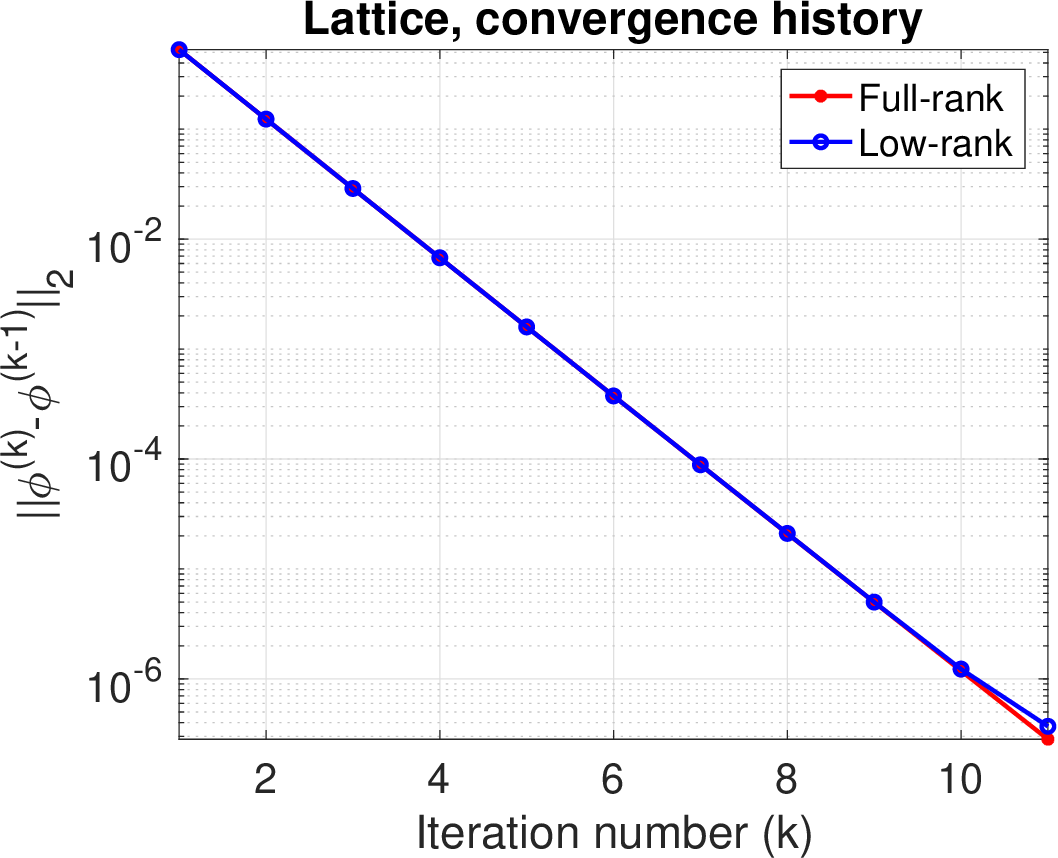}
    \includegraphics[width=0.45\textwidth]{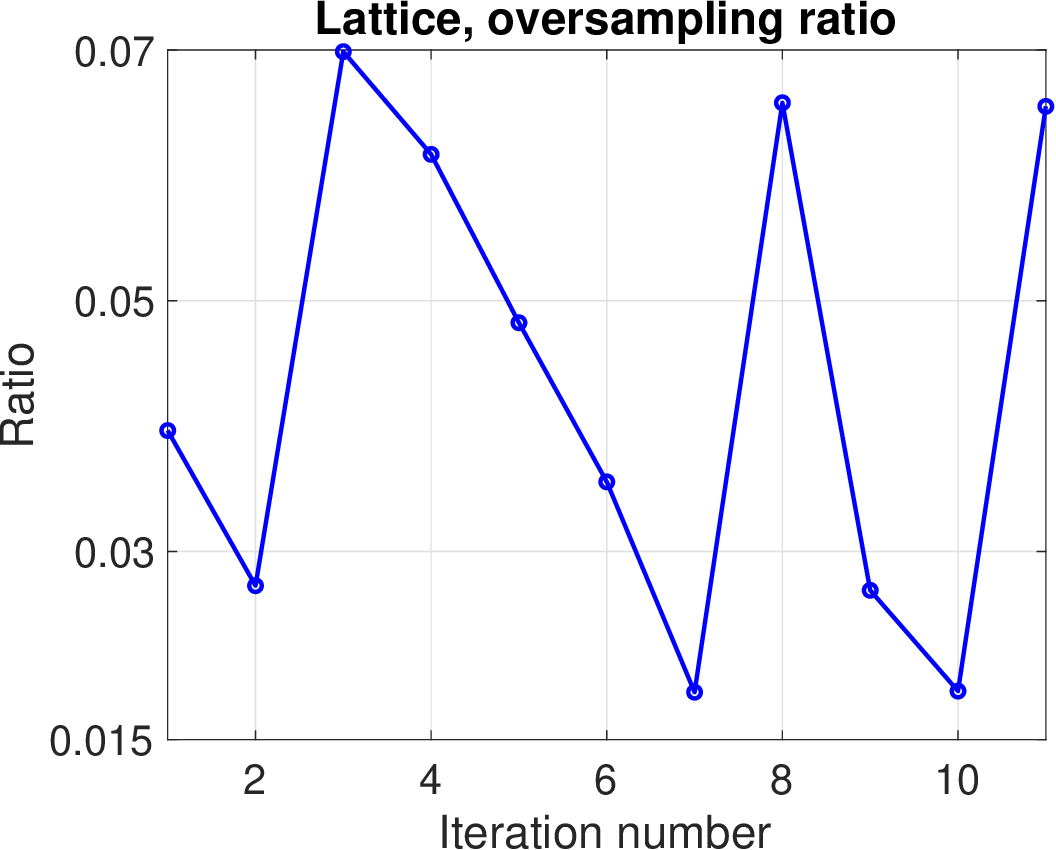}
  \caption{Convergence history and oversampling ratio during low-rank source iterations for the lattice problem in Sec. \ref{sec:lattice}. Left: convergence history.  Right: oversampling ratio. \label{fig:lattice-performance}}
\end{figure}

\pzc{\textbf{Robustness of the low-rank SI-DSA.}}
As shown in the left plot of Fig. \ref{fig:lattice-performance}, the histories of \(\|\bphi^{(k)}-\bphi^{(k-1)}\|_2\) for the low-rank \pzc{SI-DSA} and the full-rank SI-DSA almost coincide. Similar observations can be made for the convergence histories using MGS and full QR for basis update.  As reported in Tab. \ref{tab:lattice}, \review{regardless of using MGS2-A, MGS or full QR to update the spatial basis}, the proposed low-rank SI-DSA always achieves approximately $45\%$ compression in DOFs required to store the angular flux, converges in the same number of iterations as the full-rank SI-DSA, \pzc{and almost the same error with respect to the full-rank reference solution}.

\textbf{Necessity of \review{MGS2-A and} incremental update.} \review{Despite close compression ratio achieved, these three basis update methods differ significantly in their efficiency. As shown in Tab. \ref{tab:lattice}, compared to full-rank SI-DSA,  low-rank SI-DSA achieves approximately \review{$1.67\times$} speedup using MGS2-A, but only $1.15\times$ speedup using single MGS projection with global reorthogonlization, and even becomes slower when using full-QR.}  \review{This significant performance difference is due to the number of full-QR pass required in the basis update. Using full QR in the basis update of every single inner-loop iteration leads to a total number of $5521$, while MGS with global reorthognalization reduces this number to $655$ and \review{MGS2-A totally removes it}.}

\textbf{Verification of mild augmentation.} As shown in the right panel of Fig. \ref{fig:lattice}, the iterative incremental update of the spatial basis effectively avoids aggressive augmentation during the low-rank SI iterations, resulting in less than \review{$7\%$ over-augmentation in  all iterations.}

\review{\textbf{Influence of inner-loop low-rank ``truncation" thresholds.} We also test the impact of thresholds, $\epsilon_{\textrm{res.}}$ and $\epsilon_{\textrm{diff.}}$ on the performance of our method. Specifically, we fix the outer-loop SI-DSA convergence tolerance as $\epsilon_{\textrm{SI-SA}}=5\times10^{-5}$, scan $\epsilon_{\textrm{res.}}=\epsilon_{\textrm{diff.}}=10^{-3},10^{-4},\dots,10^{-7}$ and set $\epsilon_{\textrm{SVD}}=0.1\epsilon_{\textrm{res.}}$. The results are summarized in Tab. \ref{tab:lattice-truncation-scan}, and our main observations are as follows.}
\begin{enumerate}
    \item \review{The accuracy of the low-rank method, measured by the
    difference between the low-rank and full-rank solutions, is controlled
    by the thresholds $\epsilon_{\textrm{res.}}=\epsilon_{\textrm{diff.}}$.
    When these thresholds are larger than $\epsilon_{\textrm{SI-DSA}}$,
    the discrepancy between the low-rank and full-rank solutions can also
    exceed $\epsilon_{\textrm{SI-DSA}}$.}
    \item \review{Larger thresholds $\epsilon_{\textrm{res.}},\epsilon_{\textrm{diff.}}$ leads to lower effective rank and higher compression. For moderately large thresholds, higher compression further results in greater speedup. However, when the thresholds are too large, the outer iteration of the low-rank SI-DSA may fail to converge within $25$ outer-loop iterations (see  $\epsilon_{\textrm{diff.}}=10^{-3}$).}
\end{enumerate}

\begin{table}
\centering
\scriptsize
\begin{tabular}{|c|c|c|c|c|c|c|c|c|c|}
\hline
&  Rank & $\|\bpsi_{\textrm{LR}}-\bpsi_{\textrm{FR}}\|_2$ & $\|\bphi_{\textrm{LR}}-\bphi_{\textrm{FR}}\|_2$ & C-R & FR-Iter & LR-Iter & Speedup \\
\hline
\review{MGS2-A} & \review{$229$} & \review{$3.75 \times 10^{-7}$} & \review{$3.51 \times 10^{-8}$} & \review{$45.62\%$} & \review{$11$} & \review{$11$} & \review{$1.67\times$} \\
\hline
\review{MGS-QR} & $227$ & $4.63 \times 10^{-7}$ & $3.27 \times 10^{-8}$ & $45.22\%$ & $11$ & $11$ & $1.15\times$ \\ 
\hline
QR & $227$ & $4.63 \times 10^{-7}$ & $3.27 \times 10^{-8}$ & $45.22\%$ & $11$ & $11$ & $0.65\times$ \\ \hline
\end{tabular}
\\ 
\medskip
\begin{tabular}{|p{0.28\textwidth}|c|c|c|}
\hline
& \centering\review{MGS2-A} & \review{MGS-QR} & \review{QR} \\
\hline
\review{Total number of full QR in the inner-loop basis update}
& \review{$0$} & \review{$655$} & \review{$5521$} \\ 
\hline
\end{tabular}
\caption{\review{Results for the lattice problem in Sec. \ref{sec:lattice}. In the basis update step, MGS2-A, MGS-QR and full QR are applied. \review{MGS2-A: double modified Gram-Schmidt projections with adaptive reorthogoanlization in Alg. \ref{alg:mgs2-a}. MGS-QR: single modified Gram-Schimdt with full QR reorthogoanlization if necessary}. QR: direct QR decomposition. FR/LR-Iter: number of iterations required for convergence by the full/low-rank SI-DSA.} \label{tab:lattice}}
\end{table}

\begin{table}
\centering
\review{
\scriptsize
\begin{tabular}{|c|c|c|c|c|c|c|c|}
\hline
$\epsilon_{\textrm{res}.},\epsilon_{\textrm{diff}.}$ & Rank
& $\|\bpsi_{\textrm{LR}}-\bpsi_{\textrm{FR}}\|_2$
& $\|\bphi_{\textrm{LR}}-\bphi_{\textrm{FR}}\|_2$
& C-R & FR-Iter & LR-Iter & Speedup \\
\hline
$10^{-3}$ & $123$ & $8.50 \times 10^{-4}$ & $6.92 \times 10^{-5}$
& $24.50\%$ & $8$ & $25$ (\ding{55}) & $1.13\times$ \\
\hline
$10^{-4}$ & $156$ & $2.86 \times 10^{-4}$ & $1.07 \times 10^{-5}$
& $31.08\%$ & $8$ & $8$ & $2.61\times$ \\
\hline
$10^{-5}$ & $187$ & $1.01 \times 10^{-4}$ & $3.27 \times 10^{-6}$
& $37.26\%$ & $8$ & $8$ & $2.27\times$ \\
\hline
$10^{-6}$ & $213$ & $8.81 \times 10^{-6}$ & $4.21 \times 10^{-7}$
& $42.44\%$ & $8$ & $8$ & $1.81\times$ \\
\hline
$10^{-7}$ & $230$ & $3.68 \times 10^{-7}$ & $2.87 \times 10^{-8}$
& $45.82\%$ & $8$ & $8$ & $1.60\times$ \\
\hline
\end{tabular}
}
\caption{\review{Results for the lattice problem with various $\epsilon_{\textrm{res}.},\epsilon_{\textrm{diff}.}$ and $\epsilon_{\textrm{SVD}}=0.1\epsilon_{\textrm{res}.}$. The convergence tolerance of SI-DSA tolerance is fixed at $5\times 10^{-5}$. The symbol \ding{55} indicates failure of convergence within $25$ outer-loop source iterations.}}
\label{tab:lattice-truncation-scan}
\end{table}

\subsection{Multiscale variable scattering problem\label{sec:variable-scattering}}
\begin{figure}[]
\centering
    \hspace{-0.5cm}
    \includegraphics[trim=10mm 9mm 10mm 0mm, clip,width=0.33\textwidth]{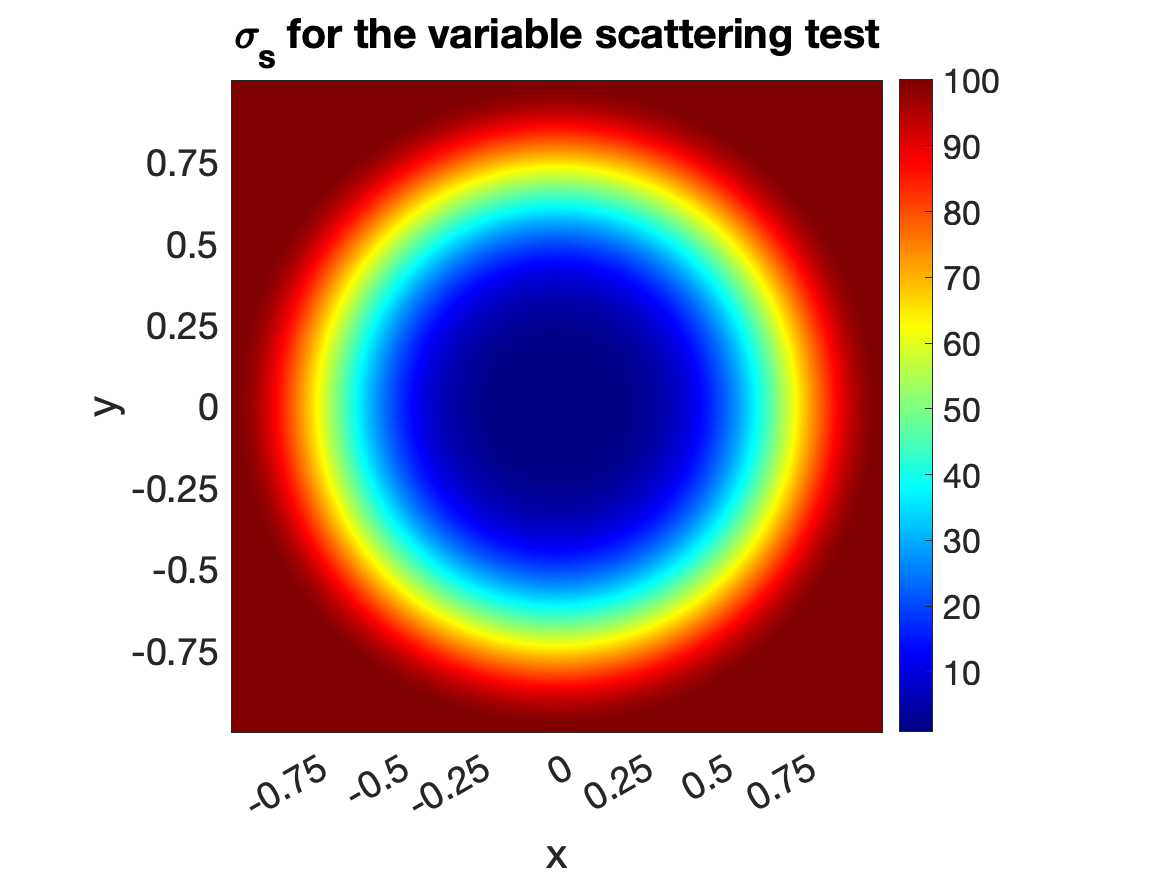}
    \hspace{-0.5cm}
    \includegraphics[trim=10mm 9mm 10mm 0mm, clip,width=0.33\textwidth]{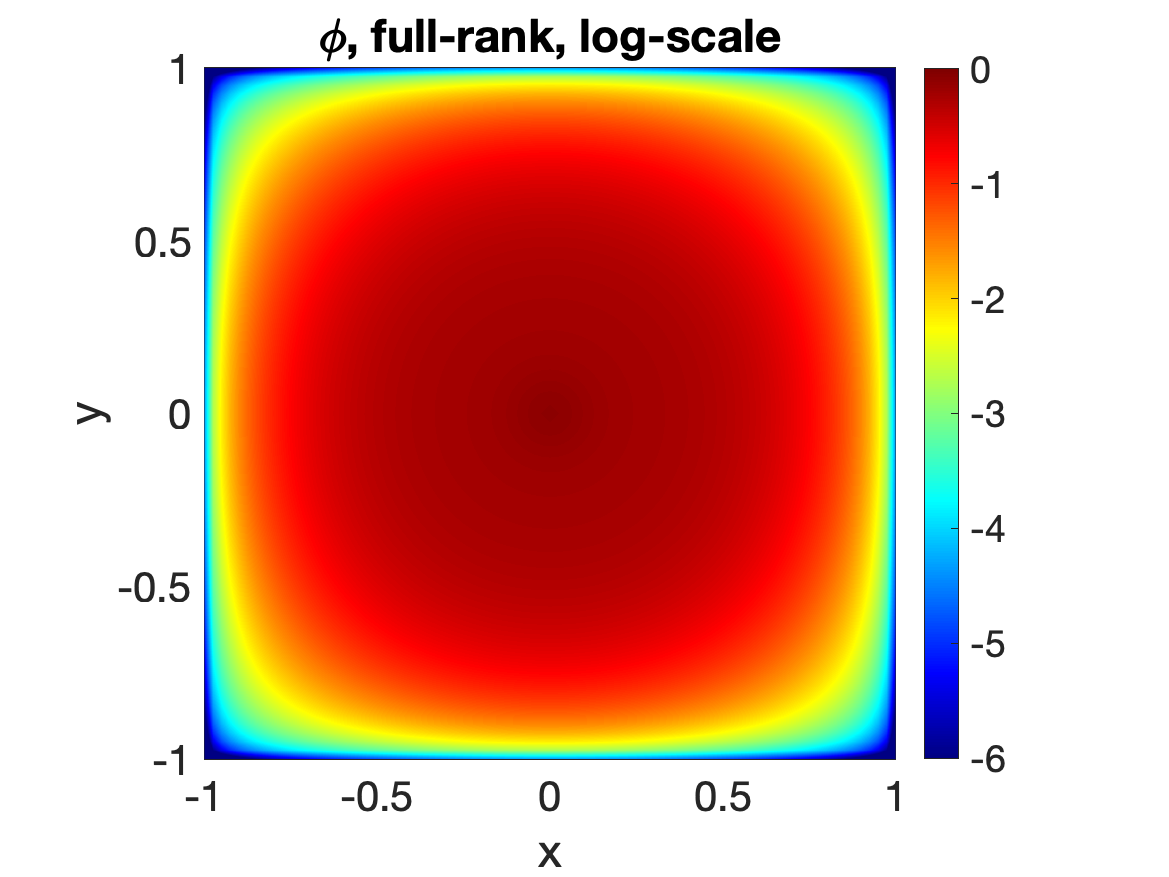}
    \hspace{-0.5cm}
    \includegraphics[trim=10mm 9mm 10mm 0mm, clip,width=0.33\textwidth]{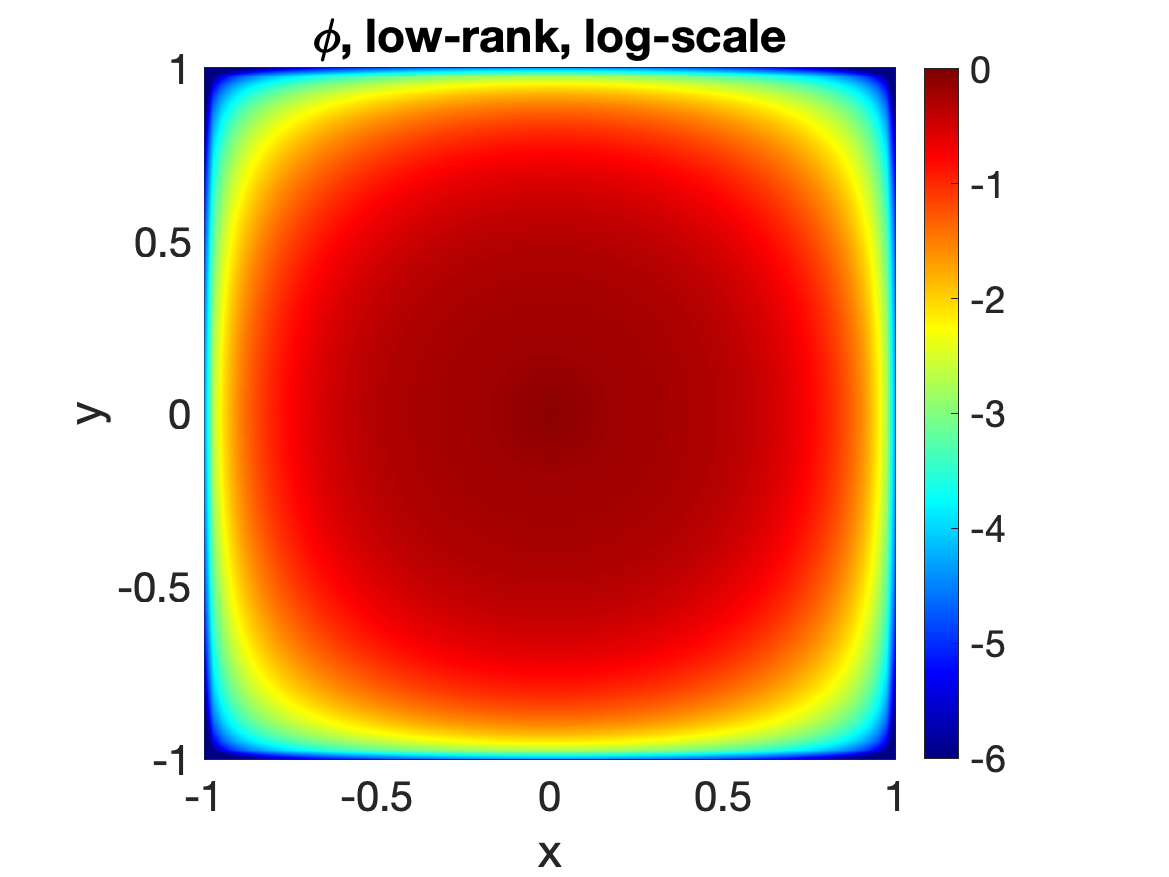}
  \caption{Configuration of scattering cross section and scalar fluxes obtained by full-rank and low-rank SI-DSA for the variable scattering example in Sec. \ref{sec:variable-scattering}. Left: $\sigma_s(x,y)$. Middle: $\phi$ obtained by full-rank SI-DSA. Right: $\phi$ obtained by low-rank SI-DSA. \label{fig:variable-scattering}}
\end{figure}

We simulate a multiscale, variable scattering problem on the computational domain $[-1,1]^2$ with vacuum boundary conditions, no absorption, and a Gaussian source $G(x,y)=\exp(-100(x^2+y^2))$. The scattering cross section (illustrated in Fig. \ref{fig:variable-scattering}) is defined as
\begin{equation}
    \sigma_s(x,y) = \begin{cases}
                    99.9r^4(r^2-2)^2+0.1, \quad\text{if}\; r=\sqrt{x^2+y^2}\leq 1,\\
                    100,\quad \text{otherwise}.
                    \end{cases}
\end{equation}
This scattering cross section increases smoothly from $0.1$ in the center to $100$ towards the boundary, implying a smooth transition from the transport regime in the center to the scattering regime near the boundary.  Hence, this example exhibits significant multiscale effects. We partition physical domain with a $80\times80$ mesh and employ the CL($40,20$) quadrature in angular space. The scalar fluxes obtained by full-rank and low-rank SI-DSA are presented on a logarithmic scale in Fig. \ref{fig:variable-scattering}. No visual difference is observed.

We further test the robustness of our low-rank SI-DSA with respect to the hyper-parameters in the angular sampling step, namely $p$ (the number of angles selected each iteration) and $q$ (the number of randomly sampled candidate angles). The results are summarized in Tab. \ref{tab:variable-scattering}. Our main observations are as follows. 

\textbf{Robust accuracy and convergence against various $p$ and $q$.} Our low-rank SI-DSA converges in the same number of iterations as the full-rank SI-DSA.  Moreover, the difference between the scalar fluxes produced by the low-rank and full-rank methods,
$\|\bphi_{\textrm{LR}}-\bphi_{\textrm{FR}}\|_2$ is on the order of $10^{-7}$-$10^{-8}$, at least one order of magnitude smaller than the threshold $\epsilon_{\textrm{SI-SA}}=10^{-6}$ used in the stopping criterion of SI-DSA, i.e. $\|\bphi^{(n)}-\bphi^{(n-1)}\|_\infty\leq \epsilon_{\textrm{SI-SA}}$ .

\begin{figure}[]
\centering
    \includegraphics[width=0.45\textwidth]{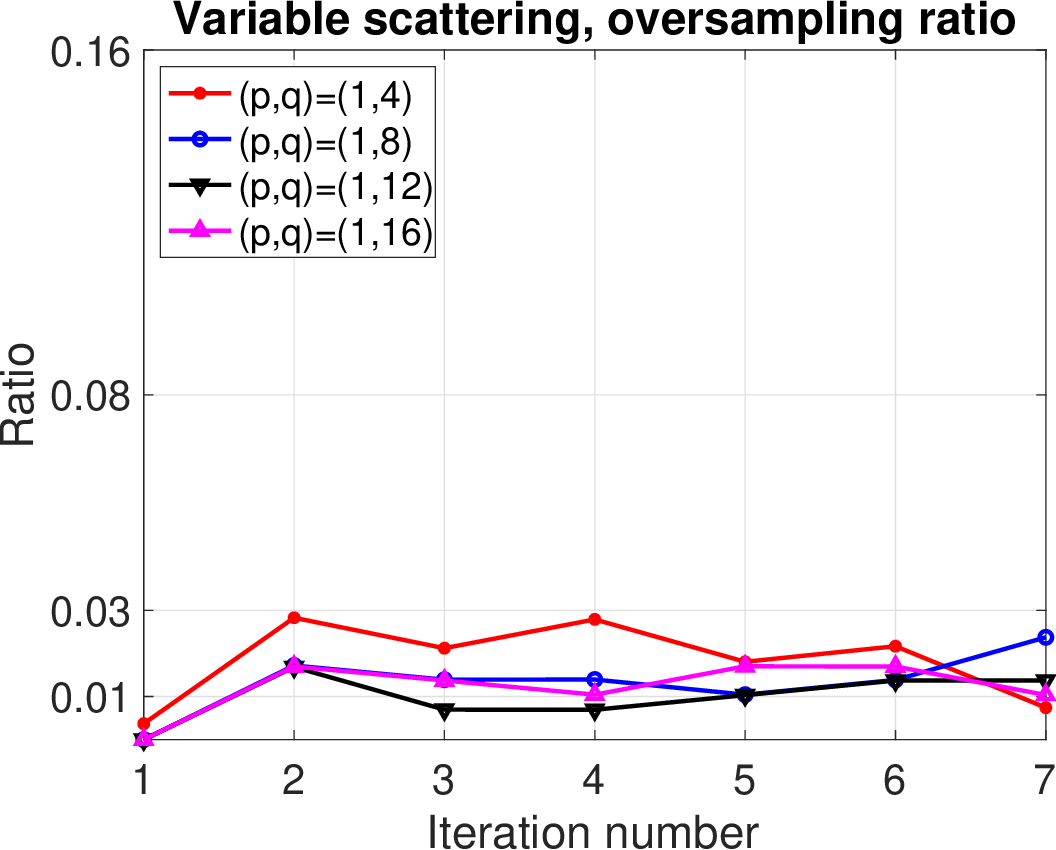}
    \includegraphics[width=0.45\textwidth]{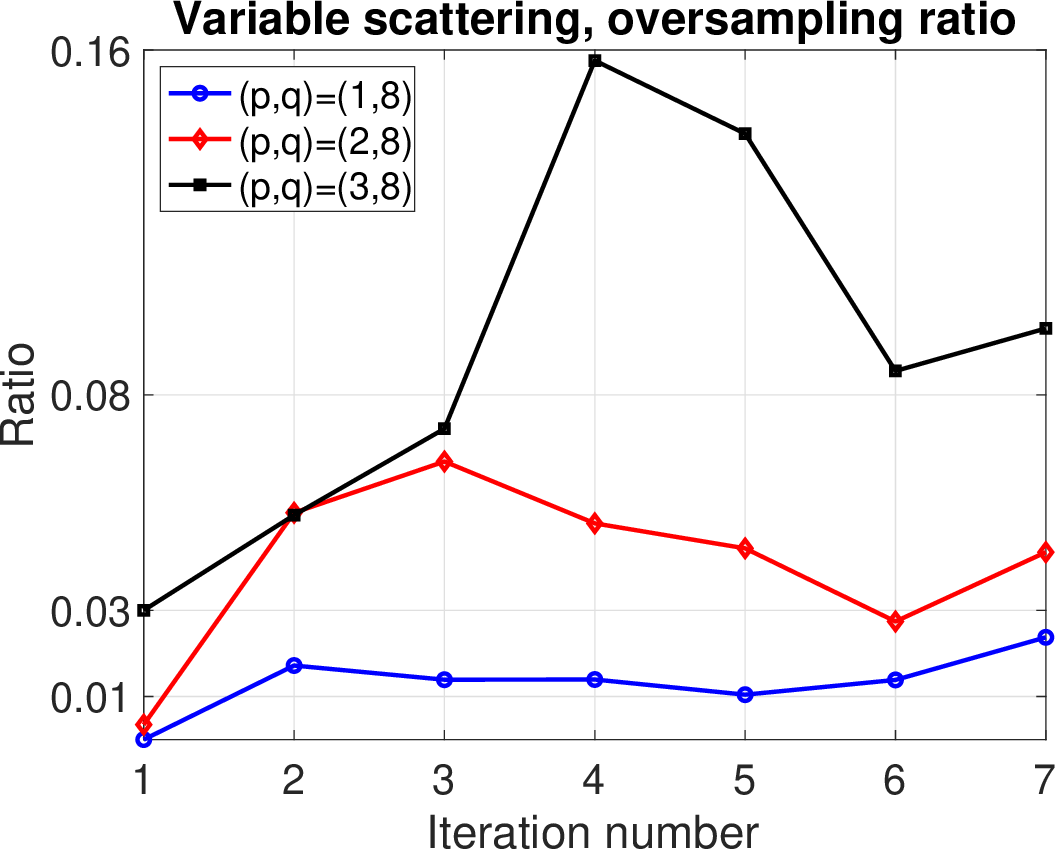}
  \caption{Oversampling ratio during low-rank SI for the variable scattering problem. Left: fixed $p=1$ and different $q$. Right: fixed $q=8$ and different $p$.\label{fig:variable-scattering-performance}}
\end{figure}

\textbf{Robust \pzc{acceleration and memory saving} against $p$ and $q$.} We observe that, for various $p$ and $q$, the low-rank SI-DSA consistently reduces the DOFs required to store the angular flux \pzc{to} approximately \pzc{$38-40\%$ of full rank}. \pzc{Moreover, for $p=1$, \review{$1.79\times-1.93\times$} speedup is achieved for $q=4,8,12,16$.  For $q=8$, speedup  is \review{$1.90\times$} for $p=1$,  \review{$1.83\times$} for $p=2$ and \review{$1.67\times$} for $p=3$, consistent with the over-sampling ratio presented in Fig. \ref{fig:variable-scattering-performance}.}

\textbf{Mild space augmentation.} In Fig. \ref{fig:variable-scattering-performance}, we present the effective rank and the oversampling ratio during low-rank SI for various $p$ and $q$. With $p=1$, the oversampling ratio remains below \review{$3\%$} for \review{$q=4,8,12,16$}. For fixed $q=8$, increasing the number of sampled angles per iteration leads to a higher oversampling ratio, but it never exceeds \review{$16\%$}. These observations verifies the effectiveness of the proposed \pzc{method} in preventing overly aggressive space augmentation during rank adaptation, \review{and implies that a small $p$ value be beneficial to control the oversampling.}  \review{We also point out that the oversampling ratio in the basis adaptation  is $100\%$ for standard augmented BUG \cite{ceruti2022rank} and even higher for standard truncated summation framework. Such oversampling may prevent practical saving in the regime of this example whose effective rank is $38\%$ of the full rank.}

\textbf{Comparison between the rank of low-rank solutions and the effective rank of the full-rank solution.} We compute the effective rank of the full-rank solution by rewriting it in the matrix format as in \eqref{eq:lr-representation} and performing truncated SVD with the same truncation tolerance as in the low-rank method, i.e., $\epsilon_{\textrm{SVD}}=10^{-8}$. The resulting effective rank is 318, which is close to the ranks of the low-rank solutions ranging from \review{270} to \review{304}.
\begin{table}[htbp]
\centering
\begin{subtable}{\textwidth}
\centering
\small
\begin{tabular}{|c|c|c|c|c|c|c|c|}
\hline
$q$ & Rank & $\|\bpsi_{\textrm{LR}}-\bpsi_{\textrm{FR}}\|_2$ & $\|\bphi_{\textrm{LR}}-\bphi_{\textrm{FR}}\|_2$ & C-R & FR-Iter & LR-Iter & Speedup \\ \hline
\review{$4$}  & \review{$270$} & \review{$4.09\times10^{-6}$} & \review{$2.98\times10^{-7}$} & \review{$34.81\%$} & \review{$7$} & \review{$7$} & \review{$1.79\times$} \\
\review{$8$}  & \review{$295$} & \review{$1.09\times10^{-7}$} & \review{$7.17\times10^{-8}$} & \review{$38.03\%$} & \review{$7$} & \review{$7$} & \review{$1.90\times$} \\
\review{$12$} & \review{$292$} & \review{$8.70\times10^{-8}$} & \review{$4.47\times10^{-8}$} & \review{$37.64\%$} & \review{$7$} & \review{$7$} & \review{$1.93\times$} \\
\review{$16$} & \review{$290$} & \review{$1.01\times10^{-7}$} & \review{$2.87\times10^{-8}$} & \review{$37.38\%$} & \review{$7$} & \review{$7$} & \review{$1.80\times$} \\ \hline
\end{tabular}
\caption{Results for $p=1$ and different $q$.}
\end{subtable}
\\
\begin{subtable}{\textwidth}
\centering
\small
\begin{tabular}{|c|c|c|c|c|c|c|c|}
\hline
$p$ & Rank & $\|\bpsi_{\textrm{LR}}-\bpsi_{\textrm{FR}}\|_2$ & $\|\bphi_{\textrm{LR}}-\bphi_{\textrm{FR}}\|_2$ & C-R & FR-Iter & LR-Iter & Speedup \\ \hline
\review{$1$} & \review{$295$} & \review{$1.09\times10^{-7}$} & \review{$7.17\times10^{-8}$} & \review{$38.03\%$} & \review{$7$} & \review{$7$} & \review{$1.90\times$} \\
\review{$2$} & \review{$299$} & \review{$9.00\times10^{-7}$} & \review{$8.79\times10^{-8}$} & \review{$38.54\%$} & \review{$7$} & \review{$7$} & \review{$1.83\times$} \\
\review{$3$} & \review{$304$} & \review{$2.63\times10^{-6}$} & \review{$9.82\times10^{-8}$} & \review{$39.19\%$} & \review{$7$} & \review{$7$} & \review{$1.67\times$} \\ \hline
\end{tabular}
\caption{Results for different $p$ and fixed $q=8$}
\end{subtable}
\caption{Results for the variable scattering example in Sec. \ref{sec:variable-scattering}. $p$: angle sampled added each inner-loop iteration. $q$: number of randomly drawn candidate samples angles. C-R: compression ratio. FR/LR-Iter: number of iterations for convergence required by full/low-rank SI-DSA.\label{tab:variable-scattering}}
\end{table}

\subsection{Pin-cell problem\label{sec:pin-cell}}
\begin{figure}[]
\centering
    \hspace{-0.5cm}
    \includegraphics[trim=10mm 9mm 10mm 0mm, clip,width=0.33\textwidth]{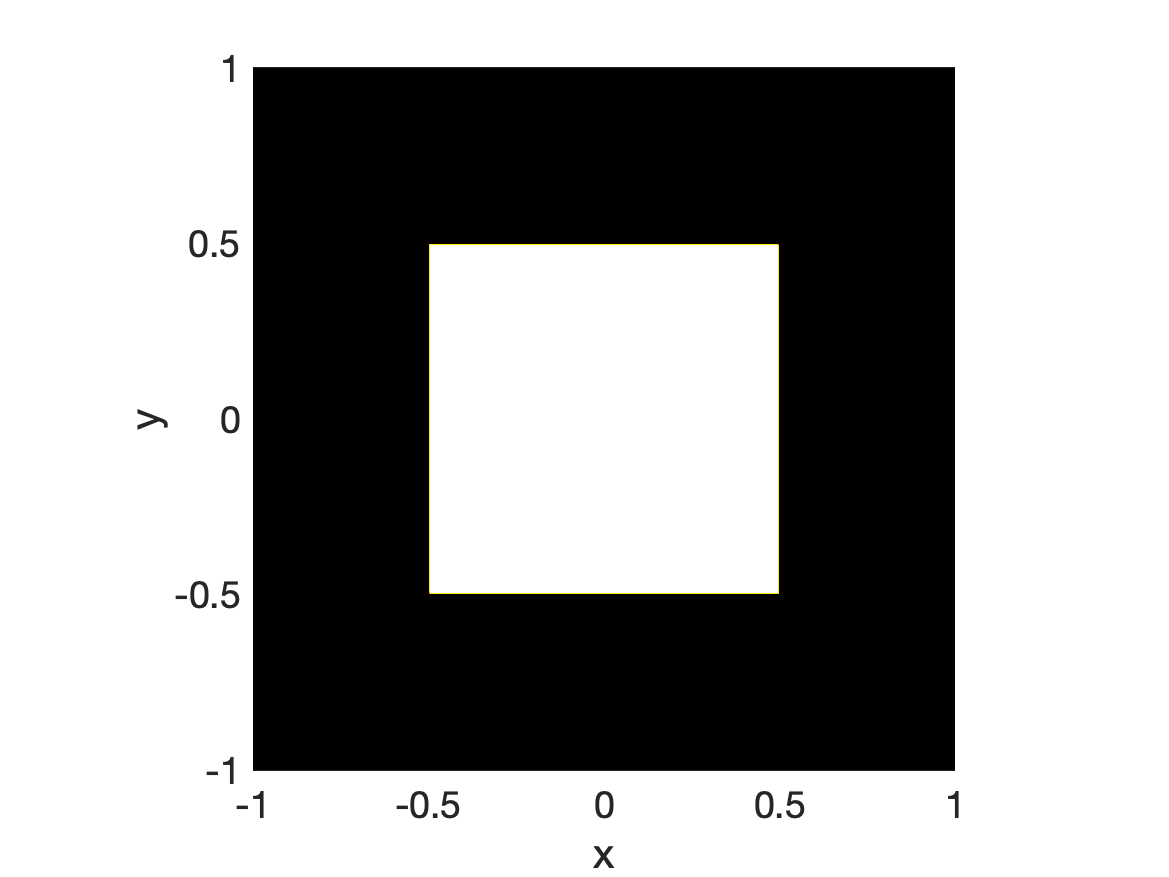}
    \hspace{-0.5cm}
    \includegraphics[trim=10mm 9mm 10mm 0mm, clip,width=0.33\textwidth]{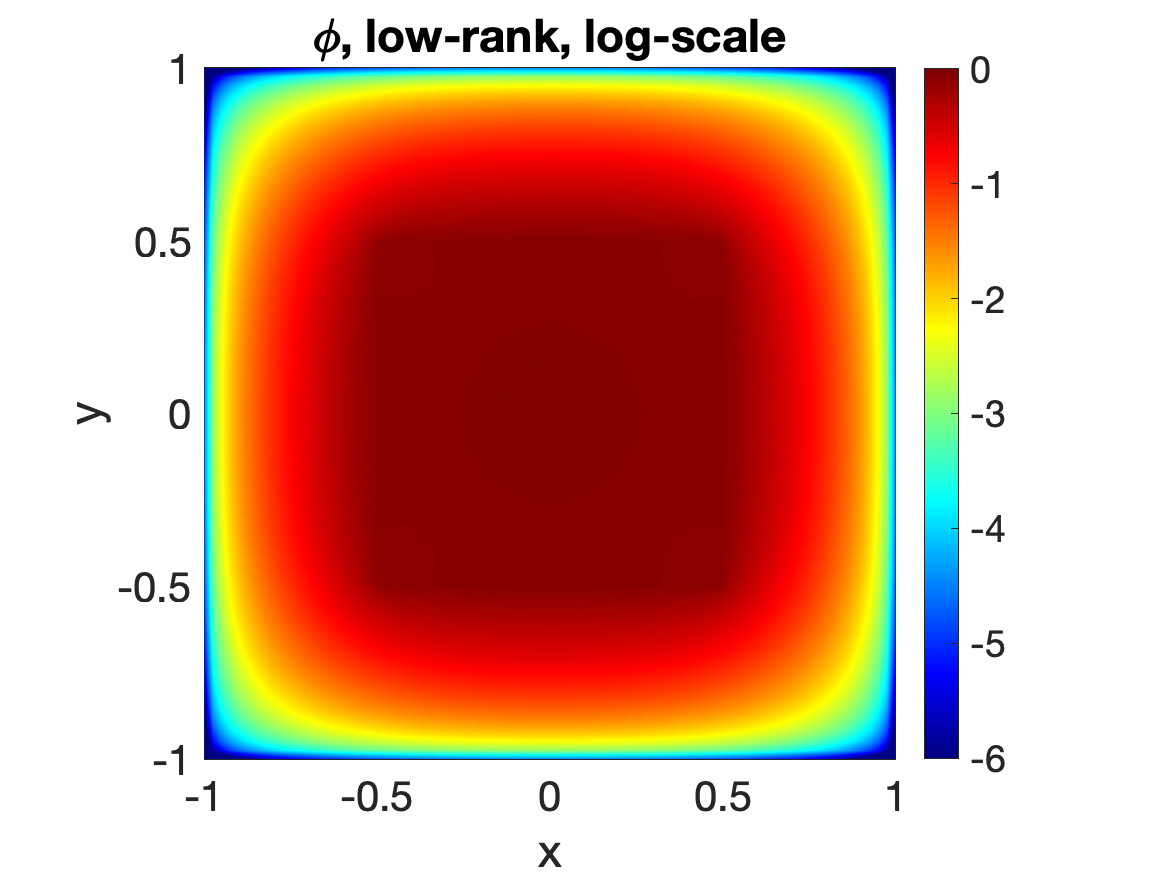}
    \hspace{-0.5cm}
    \includegraphics[trim=10mm 9mm 10mm 0mm, clip,width=0.33\textwidth]{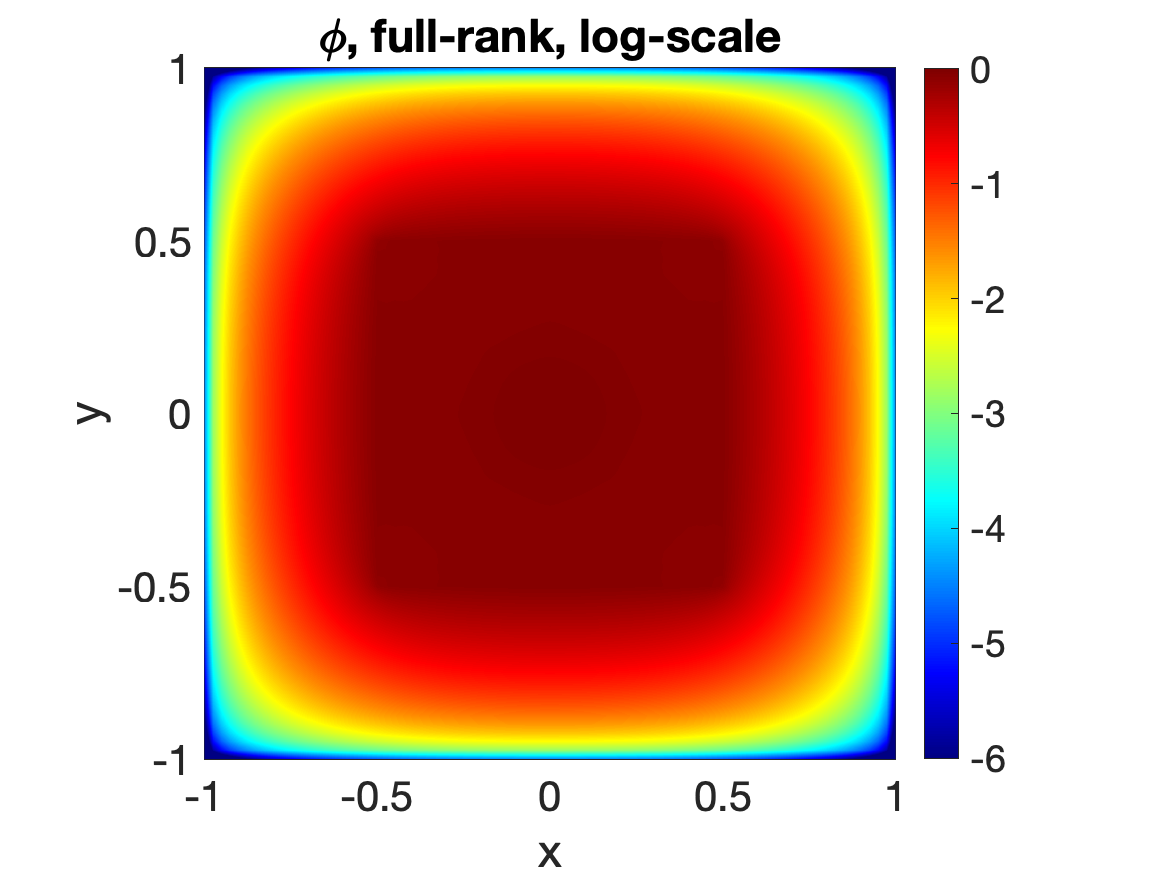}
  \caption{Configuration of scattering cross section and scalar fluxes obtained by full-rank and low-rank SI-DSA for the pin-cell example in Sec. \ref{sec:variable-scattering}. Left: configuration of $\sigma_s(x,y)$, white region corresponding to $\sigma_s=0.1$ and black region corresponding to $\sigma_s=100$. Middle: $\phi$ obtained by full-rank SI-DSA. Right: $\phi$ obtained by low-rank SI-DSA. \label{fig:pin-cell}}
\end{figure}

We consider a pin-cell problem with strong material discontinuities and multiscale effects, a classical benchmark for testing iterative solvers and preconditioners.
The computational domain is $[-1,1]^2$. There is no absorption. A Gaussian source, $G(x,y)=\exp(-100(x^2+y^2))$, and zero inflow boundary conditions are imposed. The scattering cross section (configuration illustrated in Fig. \ref{fig:pin-cell}) is defined as 
\begin{equation}
    \sigma_s(x,y)=\begin{cases}
       0.1,\; \text{if}\;|x|\leq0.5\;\text{and}\;|y|\leq0.5, \\
       100,\;\text{otherwise}.
    \end{cases}
\end{equation}
The center of the domain lies in the transport-dominated regime ($\sigma_s=0.1$), while the outer region is in the diffusion-dominated regime ($\sigma_s=100$), resulting in a 1000-fold difference in scattering strength. The domain is discretized with an $80\times80$ mesh in physical space and the CL(32,16) quadrature in angular space. Fig. \ref{fig:pin-cell} shows the scattering configuration and the scalar fluxes computed by the low-rank and full-rank SI-DSA for comparison.

The inner loop of our low-rank SI involves random subsampling of angles, with the initial angles randomly sampled in the first iteration. To assess robustness with respect to randomness, we solve the problem eight times using different random seeds determined by system time. The results are summarized in Tab. \ref{tab:pin-cell}, and the main observations are as follows.

\textbf{Effectiveness and robustness with respect to the randomness in the algorithm.} From the accuracy perspective, the difference between the solutions produced by the low-rank SI-DSA and the full-rank SI-DSA is consistently much smaller than $\epsilon_{\textrm{SI-SA}}=10^{-6}$. In \pzc{all} $8$ runs, our low-rank SI-DSA converges with exactly the same number of iterations \pzc{as its full-rank counterpart.}

The solution ranks for all runs of the low-rank SI-DSA are close to $42\%$ of the full rank. Moreover, even in this challenging scenario, thanks to the efficiency of the proposed sampling strategy and the mild space augmentation in the rank-adaptation step, the method achieves at least a \review{$2.04\times$ speedup} and an average speedup of \pzc{$2.11\times$}. In the best run, a speedup of approximately \pzc{$2.23\times$} is obtained.

\begin{figure}[]
\centering
    \includegraphics[width=0.45\textwidth]{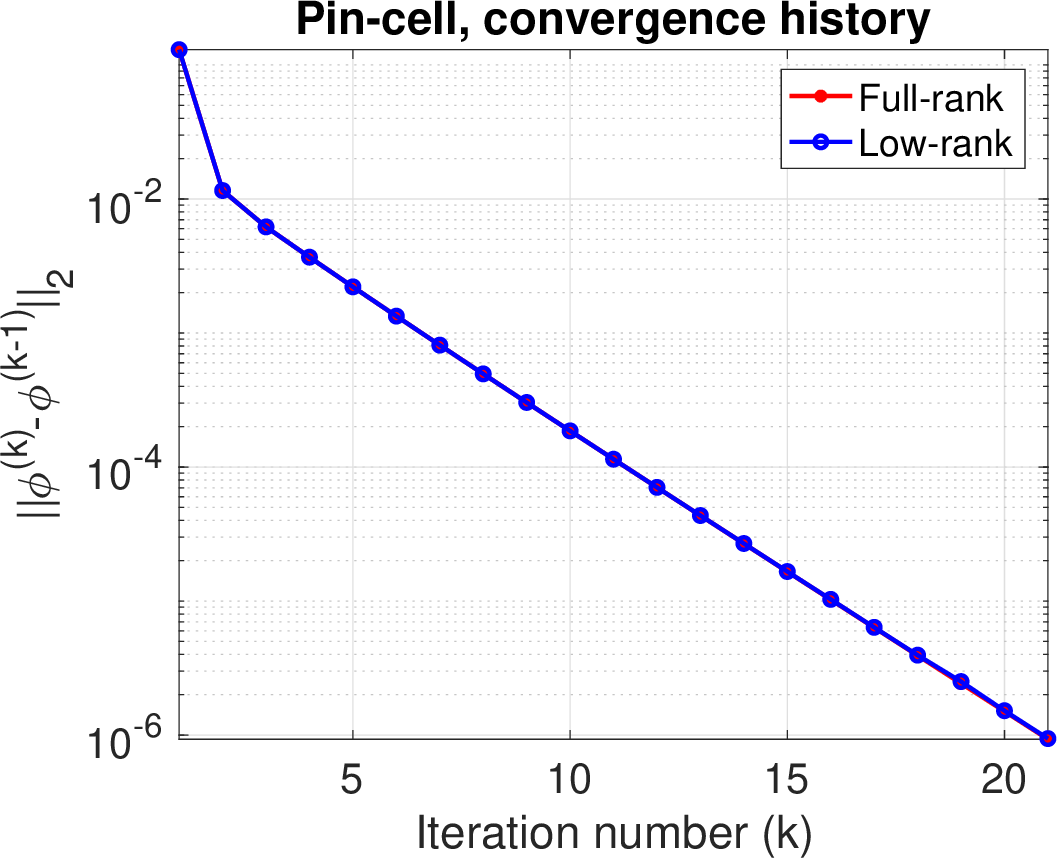}
    \includegraphics[width=0.45\textwidth]{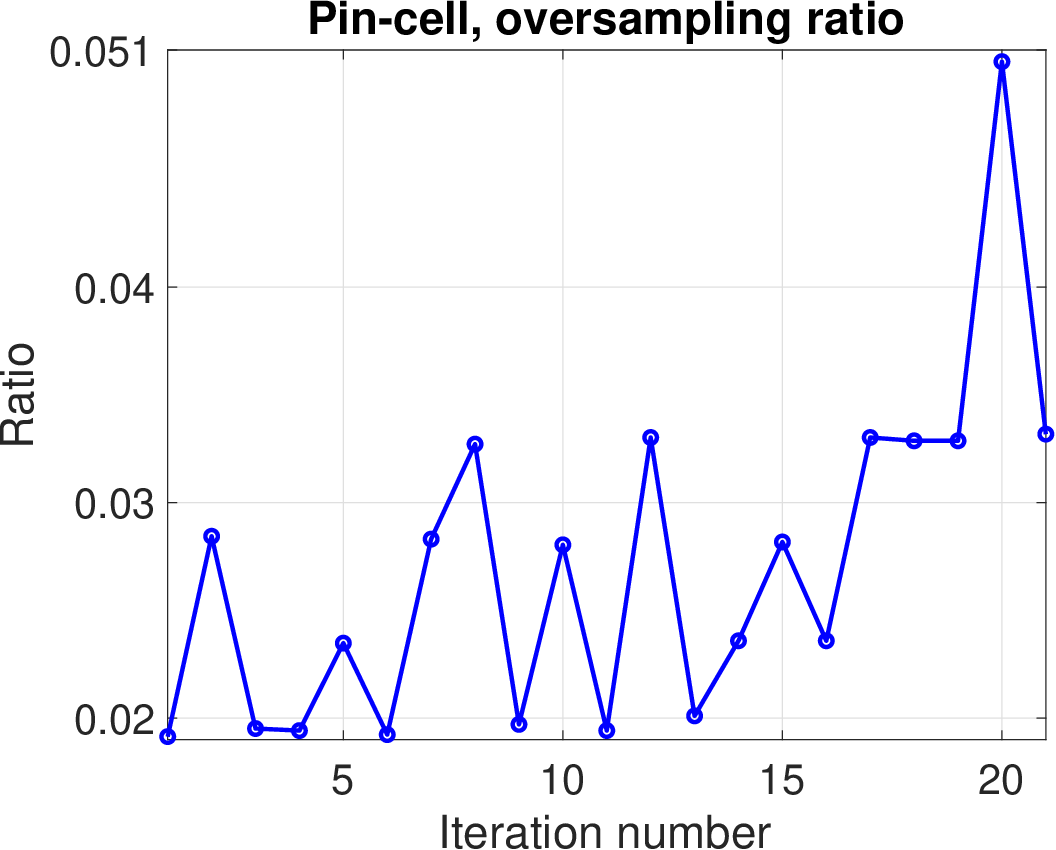}
  \caption{Convergence history and oversampling ratio during the low-rank SI for the pin-cell problem in Sec. \ref{sec:pin-cell}. Left: convergence history. Right: oversampling ratio. \label{fig:pin-cell-performance}}
\end{figure}

Figure \ref{fig:pin-cell-performance} shows the convergence history and the oversampling ratio during the low-rank SI iterations for the last run. Except for the second-to-last iteration, the change of the scalar flux in the $k$-th iteration, \(\|\bphi^{(k)}-\bphi^{(k-1)}\|_2\), for the low-rank and full-rank SI-DSA almost coincide. The oversampling ratio always remains below \review{5.1\%} across all iterations. The observations for the other runs are similar and are therefore omitted.

Finally, we compute the effective rank of the full-rank angular flux by applying truncated SVD to its matrix representation with the same truncation tolerance as in the low-rank method. The effective rank of the full-rank solution is 222, while the ranks of the low-rank solutions range from \review{$208$} to \review{$219$}, with an average of \review{$212.13$}.
\begin{table}[htbp]
\centering
\begin{tabular}{|l|c|c|c|}
\hline
Metric & Mean ($\pm$SD) & Best & Worst \\ \hline
Speedup & \review{$2.11$ ($\pm0.06$)$\times$} & \review{$2.23\times$} & \review{$2.04\times$} \\ \hline
Solution rank & \review{$212.13$ ($\pm3.87$)} & \review{$208$} & \review{$219$} \\ \hline
$\|\psi_{\text{LR}}-\psi_{\text{FR}}\|_2$ & \review{$2.99(\pm3.43)\times10^{-7}$} & \review{$8.45\times10^{-8}$} & \review{$1.12\times10^{-6}$} \\ \hline
$\|\phi_{\text{LR}}-\phi_{\text{FR}}\|_2$ & \review{$1.25(\pm0.68)\times10^{-7}$} & \review{$3.67\times10^{-8}$} & \review{$2.62\times10^{-7}$} \\ \hline
Compression ratio & \review{$42.26$ ($\pm0.77$)$\%$} & \review{$41.44\%$} & \review{$43.63\%$} \\ \hline
FR-Iter & \review{$21.00$ ($\pm0.00$)} & \review{$21$} & \review{$21$} \\ \hline
LR-Iter & \review{$21.00$ ($\pm0.00$)} & \review{$21$} & \review{$21$} \\ \hline
\end{tabular}
\caption{Statistics for the low-rank SI-DSA solving the pin-cell problem in Sec. \ref{sec:pin-cell} computed from the results of $8$ random runs. Mean: mean value. SD: standard deviation. FR/LR-Iter: number of iterations required by the convergence of the full/low-rank SI-DSA. \label{tab:pin-cell}}
\end{table}

\section{Conclusions\label{sec:conclusion}}
We develop an efficient, rank-adaptive, sweeping-based low-rank SI-DSA for solving the steady state RTE. Beyond introducing rank adaptivity into sweeping-based iterative solvers for the RTE, the main contribution lies in the design of a novel and efficient rank-adaptation strategy that updates the basis with only mild space augmentation.
Even for challenging problems with stringent outer-loop convergence criteria, i.e.,
\[
\|\bphi^{(n)}-\bphi^{(n-1)}\|_\infty \le 10^{-6},
\]
whose solution ranks are approximately $35\%–45\%$ of the full rank, the proposed method produces highly accurate low-rank solutions. The differences from the full-rank solutions are only $O(10^{-7})–O(10^{-8})$, while achieving savings in both \review{storage} and computational time. When the solution is intrinsically low-rank, e.g., in the diffusion-regime test in Sec. \ref{sec:homo}, the method achieves more than \review{$10\times$} speedup and storage. Moreover, the proposed method is shown to be robust with respect to hyperparameters and the inherent randomness in the subsampling step.

\review{
An important next step is to extend the proposed method to more complex transport applications. For time-dependent thermal radiative transfer, the proposed rank-adaptive SI-DSA solver may serve as the linear solver block within semi-implicit or fully implicit time-marching schemes, as outlined in Sec.~\ref{sec:extension}. Since the spatial variable is not further tensorized, the material–temperature equation is solved as in the classic full-rank setting, while the scalar flux and radiation energy terms are updated via contractions of the low-rank factors. For multigroup transport,  tensor train representations \cite{oseledets2011tensor} provide a natural extension to space--angle--energy decompositions, following previous work in this direction \cite{truong2023tensor,deshpande2026multigroup}. For anisotropic scattering arising in medical imaging or charged-particle transport, low-rank extensions of angular multigrid methods \cite{morel1991angular} or Fokker--Planck synthetic acceleration \cite{patel2020accelerating} will likely be needed to construct efficient preconditioners for highly forward-peaked scattering kernels.}

\review{Several important research questions also remain open. Extending both the proposed method and existing sweep-based low-rank solvers \cite{peng2023sweep,haut2026efficient} to complex geometries remains challenging, since the efficiency of the reduced operator projections in the sampling step  rely on Cartesian directional ordering. Curvilinear quadrilateral meshes may provide a feasible intermediate step toward more general geometries.  Another important direction is the development of positivity-preserving low-rank transport solvers. Recent structure-preserving low-rank methods developed for quantum dynamics \cite{appelo2025kraus} may provide useful inspiration, although the positivity constraints differ substantially from those in transport problems. Finally, conservation remains an important issue. Finite SI--DSA convergence tolerances and low-rank truncation both introduce conservation errors, making the design of efficient conservative correction techniques an important topic for future research.}

In addition, we also aim to integrate our method with more efficient or robust preconditioners, such as the quasi-diffusion method in \cite{olivier2023family} and data-driven SA preconditioners  \cite{mcclarren2022data,peng2024reduced,tang2025synthetic}.



\subsection*{Acknowledgments}
Z. Peng was partially supported by the Hong Kong Research Grants Council grants Early Career Scheme 26302724 and General Research Fund 16306825, 16307226. 
\appendix
\section{Chebyshev-Legendre quadrature\label{sec:angular-discretization}} 
The normalized CL quadrature rule, denoted by \(\mathrm{CL}(N_\theta, N_{\Omega_z})\), is the tensor product of the normalized \(N_{\Omega_z}\)-point Gauss–Legendre quadrature on [-1,1], $\{(\BOmega_{z,j},\omega_{z,j})\}_{j=1}^{N_{\Omega_z}}$ with $\sum_{j=1}^{N_{\Omega_z}}\omega_j=1$, and
the normalized $N_\theta$-point Chebyshev quadrature on the unit circle
\begin{equation}
\left\{(\theta_j,\omega_{\theta,j}): \;\theta_j = \frac{2j\pi}{N_\theta}-\frac{\pi}{N_\theta}\,\text{and}\,\omega^\phi_j=\frac{1}{N_\theta},\; j=1,\dots,N_\theta\right\}.
\end{equation}
Taking the tensor product of these two quadrature rules, we obtain the quadrature points and the corresponding quadrature weights of the CL($N_{\theta}$,$N_{\Omega_z}$) rule:
\begin{subequations}
\begin{align}
&\BOmega_j=\left(\Omega_{j,x},\Omega_{j,y},\BOmega_{j,z}\right)=\left(\cos(\theta_{j_1})\sqrt{1-\BOmega_{z,j_2}^2},\sin(\theta_{j_1})\sqrt{1-\BOmega_{z,j_2}^2},\BOmega_{z,j_2}\right),\\
&\omega_j=\omega_{\theta,j_1}\omega_{\Omega_z,j_2},
\end{align}
\end{subequations}
where $j=(j_1,j_2)$, $j_1 = 1,\ldots, N_\theta$, $j_2 =1,\ldots,N_{\Omega_z}$.
Note that if $N_\theta=2N$ and $N_{\Omega_z}=N$,  then the quadrature rule \(\mathrm{CL}(N_\theta, N_{\Omega_z})\) exactly integrates normalized polynomials of degree up to $2N$ on the unit sphere.  Moreover, CL quadrature rule is symmetric: if $\BOmega_j$ is a quadrature point, so is $-\BOmega_j$.
\section{Definitions in the matrix-vector formulation of $S_N$ upwind DG scheme \label{app:mat-vec}}
Recall that $\{\eta_i\}_{i=1}^{N_{\bx}}$ forms an orthonormal basis of $U_h^K$. For each basis $\eta_i$, we also assume that it is compactly supported in one spatial element $\mathcal{T}_{a_i,b_i}$. This definition implies that $\eta_i(x_{\half}^-,y)=\eta_i(x_{N+\half}^+,y)=\eta_i(x,y_\half^-)=\eta_i(x,y_{N+\half}^+)=0$. Following the \textit{Matlab} notation, the $(k,l)$-th elements of $\BD_{x}^\pm$, $\BSigma_t$, $\BSigma_s$ in equation \eqref{eq:mat-vec} are defined as 
{\footnotesize
\begin{subequations}
\label{eq:mat-def}
\begin{align}
\BD^{\pm}_x(k,l)&=\sum_{a=1}^{N_x}\sum_{b=1}^{N_y}\Big(-\int_{\mathcal{T}_{a,b}} \partial_x\eta_k \eta_l d\bx
-\int_{y_{b-\half}}^{y_{b+\half}}\left(\eta_k(x_{a+\half}^+,y)-\eta_k(x_{a+\half}^-,y)\right)\eta_l(x_{a+\half}^\pm,y) dy\Big),\\
\BD^{\pm}_y(k,l)&=\sum_{a=1}^{N_x}\sum_{b=1}^{N_y}\left(-\int_{\mathcal{T}_{a,b}} \partial_y\eta_k \eta_l d\bx
-\int_{x_{a-\half}}^{x_{a+\half}}\left(\eta_k(x,y_{b+\half}^+,y)-\eta_k(x,y_{b+\half}^-)\right)\eta_l(x,y_{b+\half}^\pm) dx\right),\\
\BSigma_s(k,l)&=\sum_{a=1}^{N_x}\sum_{b=1}^{N_y}\int_{\mathcal{T}_{a,b}} \sigma_s\eta_k \eta_l d\bx,\quad
\BSigma_t(k,l)=\sum_{a=1}^{N_x}\sum_{b=1}^{N_y}\int_{\mathcal{T}_{a,b}} \sigma_t\eta_k \eta_l d\bx.
\end{align}
\end{subequations}
}
The $k$-th element of vectors $\bG$ and $\bg_j^{\textrm{bc}}$ is  $       \bG(k)= \sum_{a=1}^{N_x}\sum_{b=1}^{N_y}\int_{\mathcal{T}_{a,b}} G\eta_k d\bx$ and
{\footnotesize{
    \begin{align}
   \bg_{j}^{\textrm{bc}}(k)=\sum_{b=1}^{N_y} \Big(&-\max(\Omega_{j,x},0) \int_{y_{b-\half}}^{y_{b+\half}}g(x_L,y)\eta_k(x_{\half}^+,y) dy\notag\\
    &+\min(\Omega_{j,x},0) \int_{y_{b-\half}}^{y_{b+\half}}g(x_R,y)\eta_k(x_{N_x+\half}^-,y) dy
    \Big)\notag\\
    +\sum_{a=1}^{N_x} \Big(&-\max(\Omega_{j,y},0) \int_{x_{a-\half}}^{x_{a+\half}}g(x,y_B)\eta_k(x,y_{\half}^+) dx\notag\\
    &+\min(\Omega_{j,y},0) \int_{x_{a-\half}}^{x_{a+\half}}g(x,y_T)\eta_k(x,y_{N_y+\half}^-) dx
    \Big).
    \end{align}
}}

\bibliographystyle{elsarticle-num} 
\bibliography{ref}

\end{document}